\documentclass[a4paper,10pt]{amsart}
\usepackage{epsf}  
\usepackage{amsfonts, amssymb, amsthm, amsmath}  
\usepackage{hyperref}
\usepackage{bm}
\usepackage{tikz}
\usetikzlibrary{matrix,arrows,decorations.pathmorphing}
\usepackage{tikz-cd}
\usepackage{graphicx}
\usepackage{psfrag}
\newtheorem*{correspondence formula}{Correspondence formula}
\newtheorem*{Gluing formula}{Gluing formula}
\newtheorem{thm}{Theorem}

\newtheorem{defn}[thm]{Definition}  
  
\newtheorem{claim}[thm]{Claim}  
 
\newtheorem{assumption}[thm]{Assumption} 
\numberwithin{thm}{section}

\providecommand{\totl}[1]{\ensuremath{\lceil #1\rceil }}
\providecommand{\totb}[1]{\ensuremath{\underline{ #1}}}

\DeclareMathOperator{\Aut}{Aut}

\providecommand{\rof}{{}^{\phantom{f}r}_{fg}\Omega}
\newcommand{\rhf} {{}^{\phantom{f}r}_{fg}H}
\newcommand{\rh}{{}^{r}H}
\newcommand{\ie}{\text{iE}}

\newcommand{\ee}{\text{eE}}
\newcommand{\V}{V}
\newcommand{\ex}{\bold}
\providecommand {\e}[1]{\mathfrak t^{#1}}

\newcommand{\tc}[1]{\check\rvert_{#1}}

\newcommand{\Mod}{\mathcal M}

\newcommand{\Msw}{\mathcal M^{st}}

\DeclareMathOperator{\expl}{Expl}

\newcommand{\dbar}{\bar{\partial}}

\providecommand{\et}[2]{\ensuremath{\bold T^{#1}_{#2}}}
\providecommand{\lrb}[1]{\ensuremath{\left(#1\right)}}
\providecommand{\abs}[1]{\left\lvert #1\right\rvert}

\author{Brett Parker   }
\email{brettdparker@gmail.com}  
\thanks{Funded by ARC grant  DP140100296.}
  
\title{Three dimensional tropical correspondence formula}

\begin{document}
\maketitle

\begin{abstract} A tropical curve in $\mathbb R^{3}$ contributes to Gromov-Witten invariants in all genus.
 Nevertheless, we present a simple formula for how a given tropical curve contributes to Gromov-Witten invariants when we encode these invariants in a generating function with exponents of $\lambda$ recording Euler characteristic. Our main modification from the known tropical correspondence formula for rational curves is as follows:  a trivalent vertex, which before contributed a factor of $n$ to the count of zero-genus holomorphic curves, contributes a factor of $2\sin(n\lambda/2)$. 

We  explain how to calculate relative Gromov-Witten invariants using this tropical correspondence formula, and  how to obtain the absolute Gromov-Witten and Donaldson-Thomas invariants of some $3$-dimensional toric manifolds including $\mathbb CP^{3}$. The tropical correspondence formula counting Donaldson-Thomas invariants replaces  $n$ by $i^{-(1+n)}q^{n/2}+i^{1+n}q^{-n/2}$. 
\end{abstract}

\tableofcontents

\section{Introduction}

Consider a holomorphic map of a punctured Riemann surface $\Sigma$ to $(\mathbb C^{*})^{3}$ given by $3$ meromorphic functions. To each puncture of such a surface, we can associate an integral vector $\alpha=(\alpha_{1},\alpha_{2},\alpha_{3})$, where in a holomorphic coordinate $z$ centered on the puncture, our map is $(h_{1}z^{\alpha_{1}},h_{2}z^{\alpha_{2}},h_{3}z^{\alpha_{3}})$ for some non-vanishing functions $h_{i}$. The image of such a holomorphic map under the projection $(\mathbb C^{*})^{3}\longrightarrow \mathbb R^{3}$ given by $(z_{1},z_{2},z_{3})\mapsto(-\log\abs{z_{1}},-\log\abs{z_{2}},-\log\abs{z_{3}})$ is an amoeba approximating a tropical curve with a semi-infinite end traveling in the direction of $\alpha$ for each puncture of our curve.
 
\includegraphics{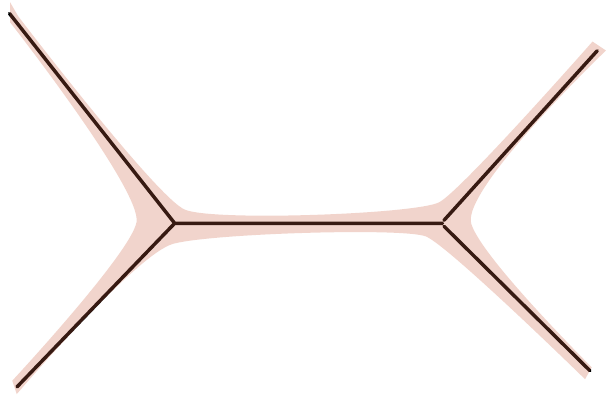}
 
 The subject of this paper is the following
 \begin{correspondence formula}[relative version] 
 \[\int_{[\mathcal M_{\bm\alpha}]}\lambda^{2g-2+n}ev_{\bm\alpha}^{*}\theta=W_{\bm\alpha}(\totb{\theta})\]
 \end{correspondence formula}
 The left hand side of this correspondence formula indicates a Gromov-Witten invariant defined as an integral over some compactification of the moduli stack of curves in $(\mathbb C^{*})^{3}$ with behavior at punctures specified by a list of integral $3$-vectors, $\bm\alpha$. The exponent of the dummy variable $\lambda$ involves the genus $g$, and number $n$ of punctures  of such curves. The right hand side indicates a weighted count of tropical curves with ends specified by this same list of integral $3$-vectors. The contribution of each tropical curve $\gamma$ to $W_{\bm\alpha}(\totb{\theta})$ is the product of a combinatorial factor determined tropically  and a weight $F_{\gamma}$. In the known tropical correspondence formula from \cite{Mikhalkin,Sieberttrop}, valid in dimension $2$, or for genus $0$, $F_{\gamma}$ may be taken as $\lambda^{2g-2+n}$, (or $1$ when not recording Euler characteristic).  We determine $F_{\gamma}$ from simpler invariants   using the following
 \begin{Gluing formula}
 \[F_{\gamma}=\sum_{\gamma'}m_{\gamma'}\prod_{v\in\V(\gamma)}\frac 1{\abs{\Aut \gamma'_{v}}}F_{\gamma'_{v}}\]
 \end{Gluing formula}
This gluing formula involves deformations $\gamma'$ of our tropical curve $\gamma$ in the case that $\gamma$ is not transverse. The combinatorial factor $m_{\gamma'}$ is determined tropically. The remaining terms involve the weights $F_{\gamma'_{v}}$ associated to simpler tropical curves $\gamma'_{v}$ obtained from each vertex $v$ of $\gamma$. We can recursively calculate $F_{\gamma}$ using this gluing formula and the following
 
 \begin{thm}\label{vertex theorem} Suppose that $\gamma$ has a single vertex and three edges labeled by integral vectors $\alpha,\beta$ and $(-\alpha-\beta)$. If $\alpha\wedge \beta$ is $n>0$ times a primitive integral vector, then
 \[F_{\gamma}=\frac{\sin (n\lambda/2)}{(n/2)}\ .\]
 If $\beta=0\neq \alpha$, then $F_{\gamma}=\lambda$.
 \end{thm}

As $(\mathbb C^{*})^{3}$ is not compact, we must specify what we mean by the Gromov-Witten invariants appearing above in our correspondence formula. There is a canonical compactification of $(\mathbb C^{*})^{3}$ to an exploded manifold\footnote{For an introduction to exploded manifolds, see \cite{iec}. For a shorter introduction in the style of lecture notes with a link to videos, see \cite{scgp}.} $\ex T^{3}$. The $[\mathcal M_{\bm\alpha}]$ in our correspondence formula indicates the virtual fundamental class of the moduli stack of curves in  $\ex T^{3}$, constructed in \cite{vfc}. Alternately, our Gromov-Witten invariants reflect the Gromov-Witten invariants of any toric compactification of $(\mathbb C^{*})^{3}$ relative to its toric boundary divisor;\footnote{ As we discuss later, suitably interpreted, these relative Gromov-Witten invariants do not depend on such a choice of compactification. We use exploded manifolds to define these Gromov-Witten invariants, however log geometry may be used as in \cite{acgw, Chen, GSlogGW}, or the geometric approach from \cite{IonelGW} may also be used.} here, the integral $3$-vector at each puncture translates to contact data with the boundary divisors. There is yet another interpretation for the invariants in Theorem \ref{vertex theorem},  using a $(\mathbb C^{*})^{3}$-equivariant construction of the virtual fundamental class of the moduli stack of such 3-punctured curves in $(\mathbb C^{*})^{3}$.

These Gromov-Witten invariants of $\ex T^{3}$ are useful in calculating absolute Gromov-Witten invariants of compact toric $3$-manifolds. For a simple relationship between our relative Gromov-Witten invariants and  absolute Gromov-Witten invariants, we restrict to toric manifolds such as $\mathbb CP^{3}$ and $(\mathbb CP^{1})^{3}$ with toric fan satisfying a convexity assumption.\footnote{Tropical correspondence formulae for Gromov-Witten invariants of other $3$-dimensional toric manifolds are complicated by contributions from curve components in the boundary, however our Gromov-Witten invariants of $\ex T^{3}$ are an essential ingredient for such general tropical correspondence formulae. Indeed, any compact $6$-dimensional symplectic manifold undergoing a suitable normal-crossing degeneration has Gromov-Witten invariants satisfying a kind of tropical correspondence formula in which our Gromov-Witten invariants of $\ex T^{3}$ play an important role.  For a  lower dimensional example, see \cite{tec}. }  The Gromov-Witten invariants of these toric manifolds enjoy a simple tropical correspondence formula, equation (\ref{absolute}) on page \pageref{absolute}. 

\

Section \ref{correspondence formula} details the  correspondence  and gluing formulae above. Worked examples of these formulae  feature in section \ref{Calculation section}, where we use these examples to determine that in the context of Theorem \ref{vertex theorem}
\[F_{\gamma}=\frac{i^{-(n+1)}q^{n/2}+i^{(n+1)}q^{-n/2}}{n}\]
where $q$ is some formal power series in $\lambda$, and $n$ is as in Theorem \ref{vertex theorem}.

%

Our tropical correspondence formula for $\ex T^{3}$ and gluing formula can not be used to determine $F_{\gamma}$ further. For that, we use known Gromov-Witten invariants of $(\mathbb CP^{1})^{3}$ and $\mathbb CP^{3}$, and a tropical correspondence formula for such Gromov-Witten invariants. Calculating some easy examples in  section \ref{calc section},  we show that $q^{\frac12}=ie^{i\lambda/2}$, so $q=-e^{i\lambda}$, which is the change of coordinates relating Gromov-Witten invariants and Donaldson-Thomas invariants,\footnote{The GW/DT correspondence formula from \cite{gwdt} relating Gromov-Witten invariants and Donaldson-Thomas invariants is the reason we express $F_{\gamma}$ in terms of $q$. We use this correspondence formula  as a black box, as a theory of Donaldson-Thomas invariants for exploded manifolds is not yet developed.} proved in the toric case \cite{gwdt}. We use this in section \ref{DT section} to give a simple tropical correspondence formula for Donaldson-Thomas invariants of some toric manifolds such as $\mathbb CP^{3}$.

  Up to the end of section \ref{DT section}, our arguments and descriptions may be understood without any background in Gromov-Witten invariants or exploded manifolds.  In section \ref{exploded section}, we explain how the tropical correspondence and gluing formulae described and explored combinatorially in  earlier sections all follow from  tropical gluing formulae for Gromov-Witten invariants of exploded manifolds. These tropical gluing formulae apply in far more general contexts than the toric  $3$-dimensional situation explored in this paper. 

 \section{Correspondence formula}
  \label{correspondence formula}

 Our correspondence formula features (balanced, parametrized, connected) tropical curves in $\mathbb R^{3}$. Such a tropical curve $\gamma$ is either an integral affine map $\mathbb R\longrightarrow \mathbb R^{3}$, or has the following data:
 \begin{itemize}
\item A finite set of vertices, $\V (\gamma)$ with a map $\V(\gamma)\longrightarrow \mathbb R^{3}$.
\item A finite set of internal edges $\ie (\gamma)$. Each internal edge $e$ is an integral-affine map $[0,l]\longrightarrow \mathbb R^{3}$, attached to a vertex  $v(e,1)\in V(\gamma)$ at $0$, and a vertex $v(e,2)$ at $l$. We require that the length $l$ of such an internal edge is strictly positive. 
\item A finite set of external edges, $\ee(\gamma)$. Each external edge $e$ is a integral-affine map $[0,\infty)\longrightarrow \mathbb R^{3}$ attached to a vertex $v(e)$ at $0$. 
 \end{itemize}
Tropical curves in this paper are balanced in the sense that that the sum of the derivatives of edges leaving a vertex is $0$.

 It is important to understand isomorphisms of tropical curves. The domain of a tropical curve is the integral-affine graph obtained by identifying  vertices with the ends of the edges they are attached to; consider the tropical curve to be an integral-affine map of this domain into $\mathbb R^{3}$. An isomorphism between such curves is an isomorphism between their domains compatible with the maps into $\mathbb R^{3}$. (In particular, the implied orientation on internal edges given by the orientation of $[0,l]$ is not part of the intrinsic data of a tropical curve.)  Throughout, we shall consider tropical curves with labelled external edges, so  $\Aut \gamma$ indicates the group of automorphisms of $\gamma$ that fix its external edges. 
 
We now construct a moduli space of tropical curves  with the same discrete data as $\gamma$. Record the position of a vertex $v$ by $x_{v}\in \mathbb R^{3}$, and the length of an internal edge $e$  by $l_{e}\in(0,\infty)$. Let   the derivative of $e$ (with a chosen orientation) be $\alpha_{e}\in\mathbb Z^{3}$.  A tropical curve with the given data  exists if  the following system of linear equations is satisfied: 
 \begin{equation}\label{equation}\alpha_{e}l_{e}-x_{v(e,1)}+x_{v(e,2)}=0 \text{ for all }e\in \ie(\gamma)\end{equation}
 Define $P_{\gamma}$ to be the solution set of the above system of equations within the space $\mathbb R^{3n}\times[0,\infty)^{k}$ that parametrizes the position of our $n$ vertices and $k$ internal edges. Each such equation is $\mathbb Z$-linear, so  $P_{\gamma}\subset \mathbb R^{3n}\times [0,\infty)^{k}$ has a natural integral affine structure.\footnote{This integral-affine space $P_{\gamma}$ only fails to be the moduli space of tropical curves with the same combinatorial type as $\gamma$ because, by describing it this way, we have labeled the vertices and edges and chosen an orientation for each edge, killing possible isomorphisms between the tropical curves parametrized by this space. The actual moduli stack is the quotient of $P_{\gamma}$ by the (integral-affine) action of some finite group. Later on, we shall restrict our attention to general tropical curves for which an evaluation map from $P_{\gamma}$ that records the position of external edges is injective. This restriction means that any isomorphism must act trivially on $P_{\gamma}$, so the quotient of  $P_{\gamma}$  by the trivial action of $\Aut \gamma$ represents the moduli stack of tropical curves with the same combinatorial type as $\gamma$.}

 \begin{defn} Say that a tropical curve is transverse if the equations (\ref{equation}) are linearly independent.
 \end{defn}
 
 For example, all zero genus tropical curves are transverse. A tropical curve containing a loop inside a plane is not transverse because the linear equations corresponding to the edges of such a loop are not linearly independent.\footnote{As we are interested in virtual counts of holomorphic curves, the condition of well-spacedness from \cite{speyer}, important for constructing holomorphic curves in the non-Archimedian setting, is not relevant to us.}

 The left-hand side of the equations (\ref{equation}) define a linear map $A:\mathbb R^{3n+k}\longrightarrow \mathbb R^{3k}$, where $k$ is the number of internal edges of $\gamma$. Our tropical curve is transverse if and only if $A$ is surjective. Our map $A$ is also $\mathbb Z$-linear, so $A(\mathbb Z^{3n+k})$ is a finite index subgroup of $\mathbb Z^{3k}$ when $\gamma$ is transverse. In this case, define
 \begin{equation}\label{mgamma}m_{\gamma}:=\abs{\mathbb Z^{3k}/A(\mathbb Z^{3n+k})}\ .\end{equation}
$m_{\gamma}=1$ if $\gamma$ has zero genus. For genus $g$ tropical curves, $m_{\gamma}$ may alternately be computed as follows.  Choose a maximal subtree of the internal edges, and use the corresponding equations (\ref{equation}) to solve for all $x_{v}-x_{v_{0}}$ in terms of these $l_{e}$. (Translational symmetry ensures that $x_{v_{0}}$ varies freely). The $g$ remaining internal edges provide $g$ ($3$-dimensional) equations in the form
 \[\sum \alpha_{e}l_{e}=0\]
 where the sum is over the edges in a loop containing one of the remaining edges.  Let $A':\mathbb Z^{k}\longrightarrow\mathbb Z^{3g}$ be the map corresponding to the left hand side of the above $g$ equations. We may equivalently calculate $m_{\gamma}$ as $\abs{\mathbb Z^{3g}/A'(\mathbb Z^{k})}$.
 
 We now define an evaluation map recording the position of an edge of a tropical curve. For any nonzero integral vector $\alpha\in\mathbb Z^{3}$, there is a projection $\pi_{\alpha}:\mathbb R^{3}\longrightarrow \mathbb R^{2}$ with kernel spanned by $\alpha$. We can identify this quotient with $\mathbb R^{2}$ so that the derivative of this projection map sends $\mathbb Z^{3}$ surjectively onto $\mathbb Z^{2}$. This is the natural integral affine structure on the quotient of $\mathbb R^{3}$.\footnote{When $\alpha$ is $k>1$ times a primitive vector, the correct quotient for some purposes is  the quotient of $\mathbb R^{2}$ by a trivial action of $\mathbb Z_{k}$. We shall not be using this quotient, and extra combinatorial factors will appear in our formula as a consequence.}   For each edge $e$, there is an evaluation map $\totb {ev}_{e}$ which measures its position (but not the location of its endpoints). If the derivative $\alpha_{e}$ on this edge is $0$, then $\totb{ev}_{e}$ simply records the image of $e$ in $\mathbb R^{3}$. When $\alpha\neq 0$, $\totb{ev}_{e}$ records the image of $e$ under the projection $\pi_{\alpha_{e}}$. Use the notation
 \[\bm\alpha(\gamma):=\prod_{e\in \ee(\gamma)}\alpha_{e}\text{ and }\totb{ev}_{\bm\alpha(\gamma)}:=\prod_{e\in \ee(\gamma)}\totb{ev}_{e}\]

 \begin{defn}\label{general def} Say that a tropical curve $\gamma$ is general if the dimension of the solution set $P_{\gamma}$ of the equations (\ref{equation}) is the number of external edges of $\gamma$, and the linear map $\totb{ev}_{\bm\alpha(\gamma)}$ is injective. 
 \end{defn}
 
 We can read off some useful facts from the above definition of a general tropical curve. A general tropical curve can not have any internal edges $e$ with $\alpha_{e}=0$, and can't have any vertices for which all attached edges have nonzero,  co-linear derivative (because the evaluation map would then have kernel). This fact and the tropical balancing condition  imply that all vertices must be at least trivalent unless the image of $\gamma$ is a single point in $\mathbb R^{3}$. A dimension count then implies that, with the exception of curves with image a single point,  all general, transverse tropical curves are trivalent. 
 
To each general tropical curve $\gamma$ with  $n$ external edges, we shall associate a generating function 
\[F_{\gamma}=\sum _{g=0}^{\infty}n_{g,\gamma}\lambda^{2g-2+n}\]
 where $n_{g,\gamma}\in\mathbb Q$ is a certain Gromov-Witten invariant counting genus $g$ curves, discussed in Section \ref{Fgamma section}. $F_{\gamma}$ only depends on the combinatorial type of $\gamma$, and is invariant under the action of integral-affine isomorphisms of $\mathbb R^{3}$ on tropical curves. Our correspondence formula for Gromov-Witten invariants involves $F_{\gamma}$ and the solution $P_{\gamma}$ to the system of equations (\ref{equation}).  
 
 If $\gamma$ is general and transverse, the following gluing formula for $F_{\gamma}$ holds,
 \begin{equation}F_{\gamma}=m_{\gamma}\prod_{v\in \V(\gamma)}F_{\gamma_{v}}\end{equation} 
 where $m_{\gamma}$ is defined in equation (\ref{mgamma}) above, and $\gamma_{v}$ is the tropical curve with a single vertex $v$, and three external edges with derivatives the same as the three edges leaving $v$ in $\gamma$.\footnote{The general, transverse tropical curves with image a single point, and $0$,  $1$ or $2$ external edges have $F_{\gamma}=0$ for dimension reasons, so we need only describe $F_{\gamma}$ for  those transverse general tropical curves that are trivalent.} The punch line of this paper is that $F_{\gamma_{v}}=(2/n)\sin (n\lambda/2)$, where $n=\abs{\alpha\wedge \beta}$ is the smallest number so that the wedge of the derivatives of two of the edges leaving $v$ is $n$ times a primitive integral vector. 
 
 \
 
 To compute $F_{\gamma}$ for tropical curves that are general but not transverse, we  need the following notion of a $\delta$-deformation of $\gamma$.
 
\begin{defn}Choose a vector $\delta_{e}\in\mathbb R^{3}$ for each internal edge of $\gamma$. A $\delta$-deformation $\gamma'$ of $\gamma$ is a length $l'_{e}\in\mathbb R$ associated to each internal edge $e\in\ie(\gamma)$, and a general tropical curve $\gamma_{v}'$ for each vertex $v\in \V(\gamma)$ along with a labeling of its external edges  by the edges of $\gamma$ leaving $v$, so that the following conditions hold.
\begin{enumerate}
\item The derivative of each external edge of $\gamma_{v}'$ equals the derivative of the corresponding edge of $\gamma$ leaving $v$.
\item For each internal edge $e\in \ie(\gamma)$,  the following equation holds:
\begin{equation}\label{delta equation} l'_{e}\alpha_{e}-x_{v'(e,1)}+x_{v'(e,2)}=\delta_{e}\end{equation}
where $\alpha_{e}$ is the derivative of $e$, and $v'(e,i)$ is the vertex of $\gamma'_{v(e,i)}$ attached to the edge corresponding to $e$. 
\end{enumerate}
 \end{defn}
 
 The left hand side of the system of equations (\ref{delta equation}) correspond to a $\mathbb Z$-linear map $A_{\gamma'}:\mathbb R^{k}\times \prod_{v}P_{\gamma'_{v}}\longrightarrow \mathbb R^{3k}$. So long as $\delta$ is chosen generically, $A_{\gamma'}$ must be surjective for there to be a solution to (\ref{delta equation}). Recall that $P_{\gamma'_{v}}$ has a $\mathbb Z$-affine structure so there is a canonical lattice $\Lambda$ of integral tangent vectors to $\mathbb R^{k}\times \prod_{v}P_{\gamma'_{v}}$. Define $m_{\gamma'}:=\abs{\mathbb Z^{3k}/A_{\gamma'}(\Lambda)}$. We can now write a gluing formula for $F_{\gamma}$ for $\gamma$ any general tropical curve:
 \begin{equation}\label{gf}F_{\gamma}=\sum_{\gamma'}m_{\gamma'}\prod_{v\in\V(\gamma)}\frac 1{\abs{\Aut \gamma'_{v}}}F_{\gamma'_{v}}\end{equation}
 where the sum is over the combinatorial types of $\delta$-deformations $\gamma'$ of $\gamma$ for any fixed generic $\delta\in\mathbb R^{3k}$. As discussed in section \ref{gf section}, this gluing formula relies on a tropical gluing formula for Gromov-Witten invariants of exploded manifolds. In general, different generic choices of $\delta$  allow different combinatorial types of deformations to contribute to the above gluing formula. The above gluing formula can  be applied recursively  to calculate each $F_{\gamma_{v}'}$. After a finite number of applications of this gluing formula, we get a combinatorial expression for $F_{\gamma}$ in terms of the invariants for tropical curves with three external edges and a single vertex.

\

Now, we  extract some Gromov-Witten invariants of $\ex T^{3}$ from these $F_{\gamma}$. Consider tropical curves $\gamma$ with $n$ labeled external edges $e$ with derivatives $\alpha_{e}$, including exactly $k$ edges for which $\alpha_{e}=0$. Recall that  $\bm\alpha:=\{\alpha_{e_{1}},\dotsc,\alpha_{e_{n}}\}$ is shorthand for this ordered collection of  $\alpha_{e}\in\mathbb Z^{3}$. The target of our evaluation map $\totb{ev}_{\bm\alpha}$ is $\mathbb R^{2n+k}$ and the expected maximal dimension of the space of such tropical curves is $n$, so to count such curves, we need a $(n+k)$-dimensional constraint. Let $\totb{ \theta}$ be a  $(n+k)$-dimensional integral-affine subspace of $\mathbb R^{2n+k}$, and let $\Lambda_{\totb\theta}\subset \mathbb Z^{2n+k}$ be the lattice of integral vectors in  $T\totb{\theta}$. We have the freedom to choose $\totb{\theta}$ generically without changing $\Lambda_{\totb\theta}$.
 For a generic such $\totb{\theta}$, consider the following weighted count of general\footnote{It is not true that only general tropical curves satisfy a generic constraint $\totb{\theta}$. Tropical curves with higher dimensional moduli spaces may still exist, however such curves never contribute to Gromov-Witten invariants, and they are not included in our weighted count of tropical curves, $W_{\bm\alpha}(\totb\theta)$.} tropical curves
\begin{equation}\label{W}W_{\bm\alpha}(\totb\theta):=\sum_{\gamma}\frac {\abs{\mathbb Z^{2n+k}/(\totb{ev}_{\bm\alpha}(\Lambda_{\gamma})\oplus\Lambda_{\totb\theta})}}{\abs{\Aut \gamma}}F_{\gamma}\end{equation}
where $\totb{ev}_{\bm\alpha}(\Lambda_{\gamma})$ means the image of the lattice of integral tangent vectors to $P_{\gamma}$ under the derivative of the integral affine map $\totb{ev}_{\bm\alpha}$, and the sum is over (isomorphism classes of connected) general tropical curves $\gamma$ so that $\bm\alpha(\gamma)=\bm\alpha$ and $\totb{ev}_{\bm\alpha}(\gamma)\in \totb{\theta}$.
 This $W_{\bm\alpha}(\totb\theta)$ is a Gromov-Witten invariant of $\ex T^{3}$, and  does not depend on the position of $\totb\theta$, only $\Lambda_{\totb\theta}$, even though different curves contribute to the above sum when we change the position of $\totb\theta$. It is also invariant under the action of integral-affine isomorphisms of $\mathbb R^{3}$ on $\bm\alpha$ and $\totb \theta$. 
 
To describe what Gromov-Witten invariant $W_{\bm\alpha}(\totb\theta)$ is, let $[\mathcal M_{\bm\alpha}]$ be the virtual fundamental class of the moduli stack of holomorphic curves in $\ex T^{3}$ with the data $\bm\alpha$, (so  these curves have tropical parts that are tropical curves whose external edges are labeled and have derivative $\alpha_{e}$). $[\mathcal M_{\bm\alpha}]$ is constructed in \cite{vfc}.  Corresponding to $\totb{ev}_{e}$, there is an evaluation map $ev_{e}$ with target $\ex T^{3}$ when $\alpha_{e}$ is $0$, and target $\ex T^{2}$ otherwise, where $\ex T^{2}$ is the quotient of $\ex T^{3}$ by the $\ex T$ action with weight $\alpha_{e}/\abs{\alpha_{e}}$. The combined evaluation map
\[ev_{\bm\alpha}:=\prod_{e}ev_{e}:[\mathcal M_{\bm\alpha}]\longrightarrow \ex T^{2n+k}\]
is closely related to our tropical evaluation map $\totb{ev}_{\bm\alpha}:P_{\gamma}\longrightarrow \mathbb R^{2n+k}$.
 Roughly speaking, the strata of the tropical part\footnote{See \cite{iec}, section 4 for a discussion of the tropical part functor.} of $[\mathcal M_{\bm\alpha}]$ with maximal dimensional image in $\mathbb R^{2n+k}=\totb{\ex T^{2n+k}}$  are the polytopes $P_{\gamma}$ weighted appropriately by $F_{\gamma}$, and the tropical part of $ev_{e}$ is $\totb{ev_{e}}$. 
 
 Corresponding to $\totb{\theta}$, there is a unique  exploded submanifold $\ex T^{n+k}\subset \ex T^{2n+k}$, which is the subgroup with weights $\Lambda_{\totb\theta}$. Taking tropical parts gives $\totb{\ex T}^{n+k}=\mathbb R^{n+k}\subset \mathbb R^{2n+k}=\totb{\ex T}^{2n+k}$, equal to  the linear subspace spanned by $\Lambda_{\totb\theta}$. Let $\theta\in\rh^{2n}(\ex T^{2n+k})$ be the Poincare dual\footnote{See \cite{dre} for the construction of the Poincare dual and the construction of refined cohomology,  $\rh^{*}$.} to this exploded submanifold. Then we have the following correspondence formula:
 \begin{equation}\label{correspondence}\int_{[\mathcal M_{\bm\alpha}]}\lambda^{2g-2+n}ev_{\bm\alpha}^{*}\theta=W_{\bm\alpha}(\totb{\theta})\end{equation} 
or equivalently,  
\[W_{\bm\alpha}(\totb{\theta})=\sum_{g}\lambda^{2g-2+n}\int_{[\mathcal M_{g,\bm\alpha}]}ev_{\bm\alpha}^{*}\theta\]
where $[M_{g,\bm\alpha}]$ indicates the component of $[\mathcal M_{\bm\alpha}]$ containing curves of genus $g$. The proof of this correspondence formula is discussed in Section \ref{correspondence section}. There exists a correspondence formula computing Gromov-Witten invariants including gravitational descendants from counts of tropical curves (generalizing the correspondence formula from \cite{markwig}), however such counts involve non-general tropical curves, and we do not explore such a formula here.

With an upgrade to the definition of $\totb\theta$ and $W_{\bm\alpha}(\totb\theta)$, our correspondence formula (\ref{correspondence}) is also valid for any $2n$-dimensional cohomology class $\theta\in\rh^{2n}(\ex T^{2n+k})$. Given any such cohomology class, there exists some compact, complex $(2n+k)$-dimensional toric manifold $X$ and and a cohomology class $\theta'\in H^{2n}(X)$ related to $\theta$ as follows: The fan of $X$ determines a subdivision of $\mathbb R^{2n+k}$, which determines a refinement\footnote{A refinement of an exploded manifold is a complete, bijective submersion; see \cite{iec}.  } of $\ex T^{2n+k}$ equal to the explosion $\expl X$ of $X$ relative  its boundary divisors.\footnote{The explosion functor is discussed in section 5 of \cite{iec}; see also \cite{elc} for a description of the explosion functor as a base change in log geometry. } We can use the smooth part map, $\expl X\longrightarrow \totl{\expl X}=X$ to pull back $\theta'$ to $\expl X$, and use the refinement map $\expl X\longrightarrow \ex T^{2n+k}$ to pull back $\theta$ to $\expl X$. The pull back of $\theta'$ and $\theta$ agree in $\rh^{2n}(\expl X)$. To each complex $n$-dimensional boundary stratum $S$ of $X$, we can associate a multiplicity 
\[m_{S}:=\int_{S}\theta'\]
and label the corresponding $(n+k)$-dimensional stratum of the fan of $X$ with $m_{S}$. Then let $\totb{\theta}$ be the closure of the union of strata with nonzero multiplicity, or more usefully, a generic translate of this in $\mathbb R^{2n+k}$. The stratum, $\totb\theta_{S}$ corresponding to $S$ is locally a $(n+k)$-dimensional integral affine subspace of $\mathbb R^{2n+k}$ with integral tangent vectors $\Lambda_{\totb\theta_{S}}$. For a generic such translate, define a weighted count of general tropical curves as follows.
\begin{equation}\label{wdef}W_{\bm\alpha}(\totb{\theta}):=\sum_{S}m_{S}\lrb{\sum_{\gamma\in \totb{ev}_{\bm\alpha}^{-1}(\totb\theta_{S})}\frac {\abs{\mathbb Z^{2n+k}/\totb{ev}_{\bm\alpha}(\Lambda_{\gamma})\oplus\Lambda_{\totb\theta_{S}}}F_{\gamma}}{\abs{\Aut \gamma}}}\end{equation}
Our correspondence formula (\ref{correspondence}) also holds for such a weighted count of curves, so $W_{\bm\alpha}(\totb\theta)$ is independent of a choice of generic translation of $\totb{\theta}$. In fact, $\totb{\theta}$ may be regarded as a tropical sub-variety of $\mathbb R^{2n+k}$, there is an analogous definition of $W_{\alpha}(\totb{\theta})$ for tropical subvarieties $\totb{\theta}$, and  the resulting weighted count of curves does not depend on the choice of deformation of $\totb\theta$.

\section{Calculation of relations between $F_{\gamma}$}
\label{Calculation section}
 
 In this section, we illustrate some examples of our tropical correspondence formula and obtain relations between $F_{\gamma}$ for tropical curves $\gamma$ with a single vertex and $3$ external edges. The balancing condition tells us that the edges of such $3$-legged $\gamma$ have derivatives $\alpha,\beta,-(\alpha+\beta)$. Such a $\gamma$ is general unless $\alpha,\beta$ and $(-\alpha-\beta)$ are nonzero and colinear. Our first examples will prove the following 
 \begin{claim}\label{F independence} For $3$-legged $\gamma$ with derivatives $\alpha$, $\beta$ and $(-\alpha-\beta)$,  $F_{\gamma}$  only depends on $\abs{\alpha\wedge\beta}$.
 \end{claim}

 Until we have proved Claim \ref{F independence} use the notation $F(\alpha,\beta):=F_{\gamma}$. As $\gamma$ has no chosen ordering of its edges, 
\begin{equation}F(\alpha,\beta)=F(\beta,\alpha)=F(\alpha,-\alpha-\beta)\ .\end{equation}
The invariance of $F_{\gamma}$ under the action of integral affine isomorphisms of $\mathbb R^{3}$ (which follows from the definition of $F_{\gamma}$  in Section \ref{Fgamma section}), gives that for any matrix $A\in GL_{3}(\mathbb Z)$, 
\begin{equation}F(\alpha,\beta)=F(A\alpha,A\beta)\ .\end{equation}
Computing $F(\alpha,\beta)$ therefore reduces to the case that $\alpha=(k,0,0)$, $\beta=(a,b,0)$.

\

Throughout this section, we shall use the following convention for describing a $\totb{\theta}$ for defining a weighted count of curves $W_{\bm\alpha}(\totb\theta)$.  Integral affine constraints on external edges describe an integral affine subspace $\totb\theta_{0}$ of the target $\mathbb R^{N}$ of our evaluation map $\totb{ev}_{\bm\alpha}$. Let $m$ be the product of $\abs{\alpha_{e}}$ for each edge $e$ with a nontrivial constraint so that $\alpha_{e}\neq 0$. Use $\totb \theta:=m\totb \theta_{0}$ for these constraints, so the weighted count of curves with these constraints is defined as $W_{\bm\alpha}(\totb\theta)=mW_{\bm\alpha}(\totb\theta_{0})$ 

\

 For our first example, consider tropical curves with $4$ external edges having derivative $\alpha,\beta,-\alpha,-\beta$, and assume $\alpha\wedge\beta\neq 0$. There are only two types of such curves that are general and admit $\delta$-deformations for generic $\delta$,\footnote{When $\alpha+\beta$ or $\alpha-\beta$ are not primitive integral vectors, there are also general tropical curves with the same image, but with multiple internal edges. These will not admit $\delta$-deformations for generic $\delta$, so they do not contribute to Gromov-Witten invariants. Similarly, any tropical curve with a loop and all incident edges contained in a plane will not admit a generic $\delta$-deformation.} pictured below.
 
 \includegraphics{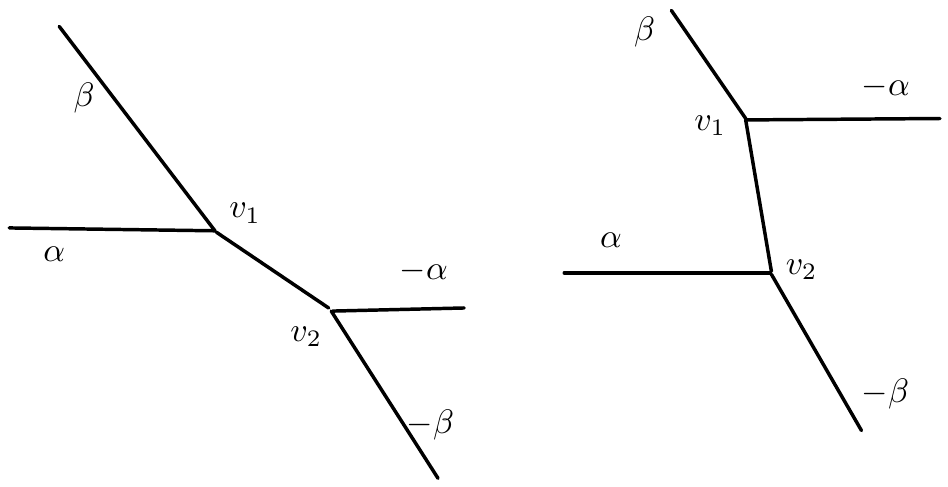}
  
  $P_{\gamma}$ has coordinates $x_{v_{1}} \in\mathbb R^{3}$ and $l_{e}\in(0,\infty)$, and $x_{v_{2}}-x_{v_{1}}$ is $l_{e}(\alpha+\beta)$ or $l_{e}(\beta-\alpha)$ in the respective cases pictured above. Our gluing formula, (\ref{gf}) implies that these tropical  curves have $F_{\gamma}$ equal to $F(\alpha,\beta)F(-\alpha,-\beta)$ and $F(-\alpha,\beta)F(\alpha,-\beta)$ respectively. Consider the weighted count of tropical curves with  the edge $\beta$ completely constrained, and  the edges $\alpha$ and $-\alpha$ constrained to lie in different planes with prescribed second coordinate.
  Constraining these edges in the two different ways pictured above give two different   calculations of the same Gromov-Witten invariant:
\[\abs{\alpha\wedge \beta}^{2}F(\alpha,\beta)F(-\alpha,-\beta)=\abs{\alpha\wedge\beta}^{2}F(-\alpha,\beta)F(\alpha,-\beta)\]
Using $GL_{3}(\mathbb Z)$ symmetry, $F(-\alpha,-\beta)=F(\alpha,\beta)$, so we obtain that for $\alpha\wedge\beta\neq 0$,
\[F(\alpha,\beta)^{2}=F(\alpha,-\beta)^{2}\] 
We show in Section \ref{zero genus section}, that $F(\alpha,\beta)=\lambda+\dotsb$, so we can remove the squares from the above equation and obtain
\[F(\alpha,\beta)=F(\alpha,-\beta)\ .\]
Using our $S_{3}$ symmetry on $F(\alpha,-\beta)$ then gives
\[F(\alpha,\beta)=F(\alpha,\beta-\alpha)\ .\]
In summary, so long as $\alpha\wedge\beta\neq 0$, we can do elementary  row operations to the matrix with rows $\alpha$, $\beta$ without affecting $F(\alpha,\beta)$. We can swap rows, add an integer multiple of one row to another, multiply a row by $-1$, and we can do the same operations on columns using $GL_{3}(\mathbb Z)$ symmetry. Using such operations, we can convert any integer matrix whose entries have greatest common divisor $k$ to a matrix in the form
\[\left[\begin{array}{ccc}k &0&0 \\ 0&nk&0\end{array}\right]\] 
where $\abs{\alpha\wedge\beta}=k^{2}n$. So, whenever $\alpha\wedge\beta\neq 0$,  $F(\alpha,\beta)$ can only depend on $\abs{\alpha\wedge\beta}$ and the greatest common divisor of $\abs{\alpha}$ and $\abs{\beta}$. 

Now consider tropical curves whose $4$ external edges have derivatives $(k,0,0)$,  $(0,nk,0)$, $(0,1,0)$, $(-k,-nk-1,0)$, where $k>0$ and $n>0 $. Again,  only two general types of such tropical curves  contribute to Gromov-Witten invariants. 

\includegraphics{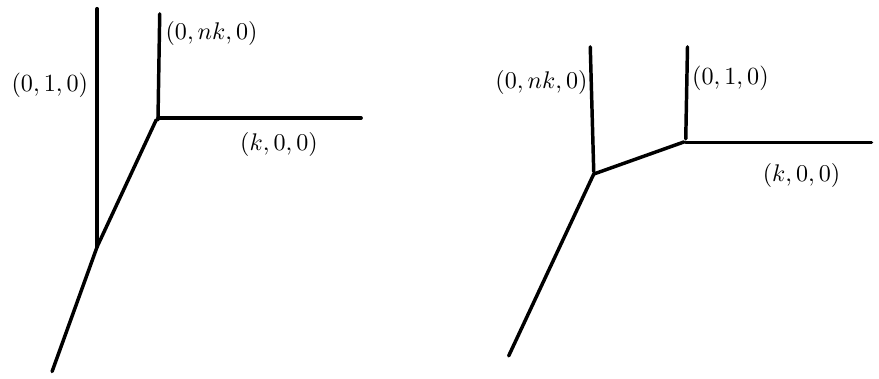}

Compute a weighted count of tropical curves with $(k,0,0)$ completely constrained and the two upward pointing edges having first coordinates constrained in the two different ways shown above. These two computations give
\[\begin{split}nk^{3}F((k,0,0),(0,nk,0))F((0,1,0),(-k,-nk-1,0))
\\ =nk^{3}F((k,0,0),(0,1,0))F((0,nk,0),(-k,-nk-1,0))\end{split}\]
however we have already determined that \[F((0,1,0),(-k,-nk-1,0))=F((k,0,0),(0,1,0))\ ,\] so 
\[F((k,0,0),(0,nk,0))=F((0,nk,0),(-k,-nk-1,0))\]
and we have proved that so long as $\alpha\wedge\beta\neq 0$,  $F(\alpha,\beta)$ depends only on $\abs{\alpha\wedge\beta}$, and  not on the greatest common divisor of $\abs \alpha$ and $\abs \beta$.  We shall explicitly calculate $F(\alpha,0)=\lambda$ in Section \ref{zero genus section}, completing the proof of Claim \ref{F independence}. 

\

We will simplify example computations that follow using the following definition of $[n]_\lambda$, which by Claim \ref{F independence}, is independent of the choice of $\alpha$ and $\beta$ so that $\abs{\alpha\wedge\beta}=n$.
\begin{equation}[n]_\lambda:=\abs{\alpha\wedge\beta}F(\alpha,\beta)\text{ where }\abs{\alpha\wedge\beta}=n\end{equation}
 Use the notation $[\alpha\wedge\beta]_{\lambda}$ for $[\abs{\alpha\wedge\beta}]_{\lambda}$. 
 There are many relations we can calculate between the different $[n]_\lambda$ from our correspondence formula. For example,  consider tropical curves with $4$ ends with derivatives $\alpha,\beta,\gamma,-\alpha-\beta-\gamma$ all in the same plane so that $\alpha\wedge\beta$, $\alpha\wedge \gamma$, $\beta\wedge \gamma$, and $(\alpha+\gamma)\wedge \beta$ all are positive rational multiples of one another.
 
 \includegraphics{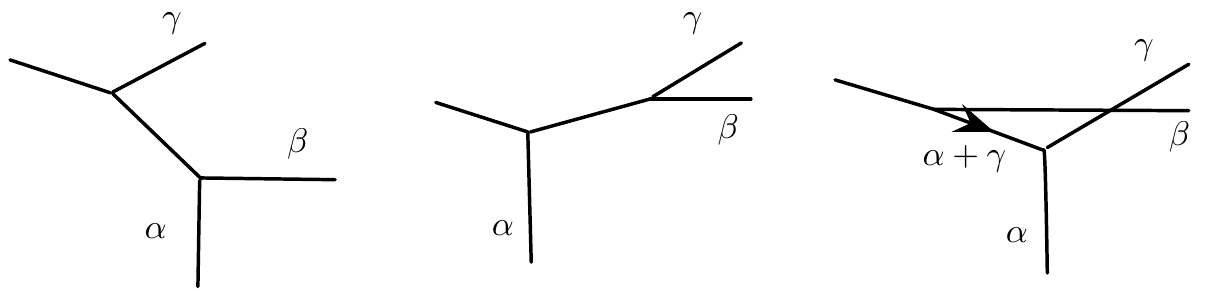}
 Equating the weighted count of tropical curves obtained constraining $\alpha$, $\beta$ and $\gamma$ as pictured to the left with the count with these edges constrained as in the other two  picture gives the following identity.
 \[[\alpha\wedge\beta]_\lambda[(\alpha+\beta)\wedge\gamma]_\lambda=[\beta\wedge\gamma]_\lambda[(\beta+\gamma)\wedge\alpha]_\lambda+[\alpha\wedge\gamma]_\lambda[(\alpha+\gamma)\wedge\beta]_\lambda\]

 \ 
 
Our next example features some non-transverse tropical curves, and leads to a beautiful formula relating the different $[n]_{\lambda}$. Consider tropical curves with $4$ ends with derivative $(1,0,0)$, $(0,1,0)$, $(-1,0,n)$, $(0,-1,-n)$. Constrain the edge in direction $(1,0,0)$  to be $(*,0,-1)$ and the   edge in direction $(0,1,0)$ to be $(0,*,1)$, then,  the only general tropical curve is shown below on the left.
 
 \includegraphics{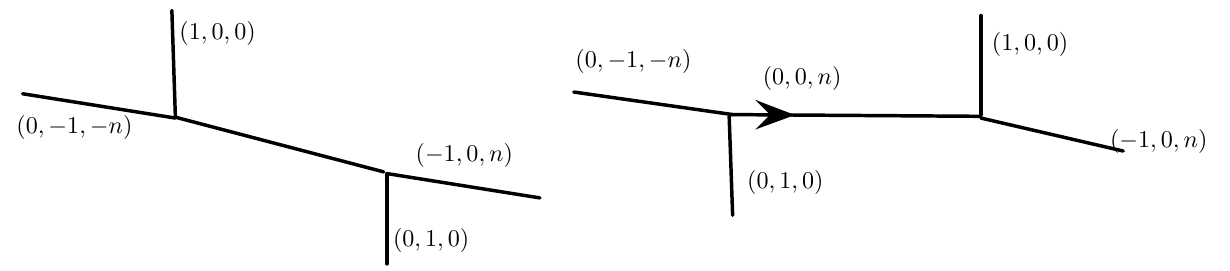}
 
  Swapping the third coordinates of the constraints on these edges allows  the general tropical curve on the right, but there are also several general tropical curves with the same image, but  several internal edges with derivative $(0,0,\mu_{i})$,
 where $\mu_{i}$ are the parts of some partition $\mu$ of $n$. Label such a tropical curve $\gamma_{\mu}$. The automorphism group of the partition $\mu$ is the same as $\Aut\gamma_{\mu}$. To calculate $F_{\gamma_{\mu}}$ using our gluing formula, we need to chose some $\delta_{e}\in\mathbb R^{3}$ for each internal edge and find $\delta$-deformations of $\gamma_{\mu}$. Using $\delta_{e_{i}}=(i,i,0)$ allows only one such $\delta$-deformation,  pictured below in the case of a partition of length $3$. The horizontal edges below look aligned because we have drawn our picture using  a projection with kernel spanned by $(1,1,0)$.
 
 \noindent \includegraphics{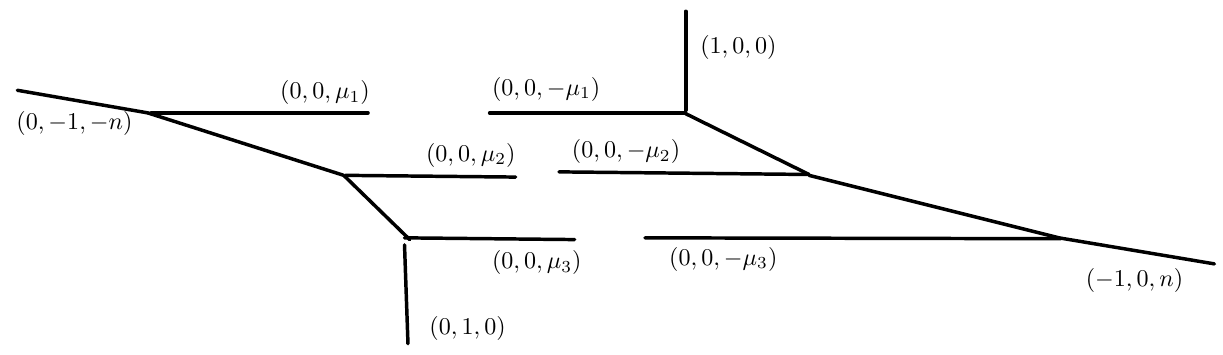}
 
For this $\delta$-deformation $\gamma'$,  $m_{\gamma'}=\prod \mu_{i}/\text{LCM}(\mu)$ where $\text{LCM}(\mu)$ is the least common multiple of $\{\mu_{i}\}$, so
  \[F_{\gamma_{\mu}}=\frac 1{LCM(\mu)} \prod_{i}\frac{[\mu_{i}]^{2}}{\mu_{i}}\]
In each case, $P_{\gamma_{\mu}}$ is $\mathbb R^{3}\times(0,\infty)$. For $(p,l)\in P_{\gamma_{\mu}}$, the position of the vertex attached to $(0,1,0)$ and $(0,-1,-n)$ is $x$, but the position of the other vertex is $x+(0,0,\text{LCM}(\mu)l)$. Equating the  weighted count of tropical curves constraining  the edges $(1,0,0)$ and $(0,1,0)$ in two different positions gives the following formula.
 \begin{equation}\label{lf}n[1]_\lambda^{2}=\sum_{\abs\mu=n}\frac 1{\abs {\Aut\mu}\prod_{i}\mu_{i}}\prod_{i}[\mu_{i}]_{\lambda}^{2}\end{equation}
 
 Some readers may recognize $\abs{\Aut \mu}\prod\mu_{i}$ as the size of the subgroup of $S_{n}$ fixing an element of cycle type $\mu$.
As the above formula is inhomogeneous,   $[n]_\lambda$ can't just be $n\lambda$; there must be higher genus contributions to $[n]_\lambda$. With sufficient combinatorial fortitude, the equations (\ref{lf}) for different $n$ can be used to derive a formula for $[n]_\lambda^{2}$ in terms of $[1]_\lambda^{2}$.  
  
 \
 
 Instead of using the above formula (\ref{lf}), we shall calculate $[n]_{\lambda}$ in terms of $[1]_{\lambda}$ by counting tropical curves with $4$ ends as pictured below. 
 
 \includegraphics{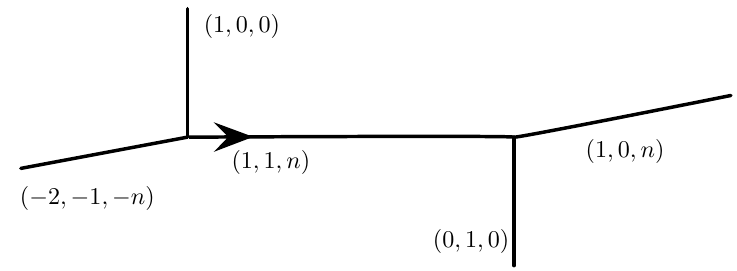}
 
 If we constrain the edges $(1,0,0)$ and $(0,1,0)$ as pictured above, there is only one such tropical curve that is general. The weighted count of this curve is $n[1]_\lambda^{2}$. On the other hand, if we constrain these edges as shown below, there are more interesting general tropical curves that contribute to our count.
 
 \includegraphics{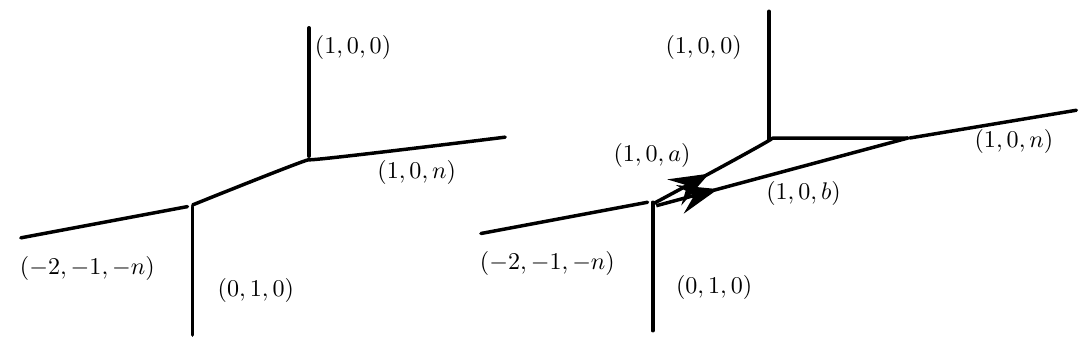}
 
 The curve shown on the right exists for all $b>a>0$ so that $a+b=n$. The tropical curves shown are the only general curves that count with nonzero weight with one exception: the derivative on the internal edge on the left-hand picture is $(2,0,n)$, twice  a primitive integral vector when $n$ is even, so when $n$ is even we may obtain a tropical curve replacing this edge with two edges, each with derivative $(1,0,n/2)$. This tropical curve could also be regarded as a degenerate case of the curve shown on the right when $a=b$. The curve on the left contributes $[1]_\lambda[n]_\lambda$ to our count if $n$ is odd, and $[2]_\lambda[n]_\lambda/2$ to our count if $n$ is even.
 
 Let us compute $F_{\gamma}$ for the curve shown above on the right. Choose $\delta_{e}$ to be $0$ for all internal edges apart from $(1,0,b)$, and let $\delta_{e}$ for that edge be $(0,1,0)$. With this suitably generic choice of $\delta$, the $\delta$-deformations of the above curve are all as pictured below.
 
 \includegraphics{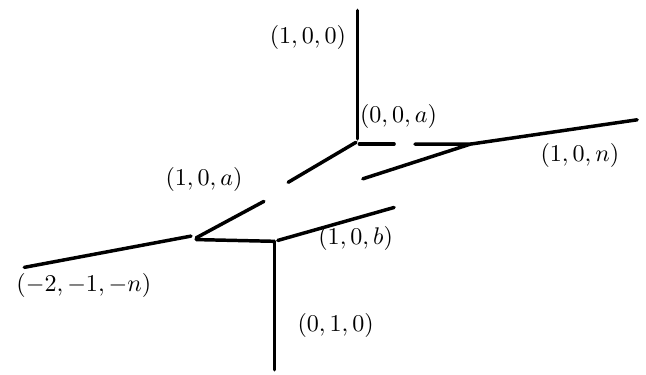}
 
 Our gluing formula then gives that $F_{\gamma}=[1]_\lambda^{2}[a]_\lambda^{2}/a$, and such a curve counts with weight $[1]_\lambda^{2}[a]_\lambda^{2}$ when we constrain the edges $(0,1,0)$ and $(1,0,0)$. Equating our two weighted counts of tropical curves in the case that $n=2m+1$ gives
 \begin{equation}\label{2m+1}(2m+1)[1]_\lambda^{2}=[1]_\lambda[2m+1]_\lambda+\sum_{a=1}^{m}[1]_\lambda^{2}[a]_\lambda^{2} \end{equation}
 
 For $n=2m$, we need to calculate the contribution of the remaining tropical curve with two internal edges of derivative $(1,0,m)$. Choose $\delta_{e}=(1,1,0)$ for the first edge, and $0$ for the second edge. Then the only $\delta$-deformed tropical curve is as pictured below. For this curve, $F_{\gamma}=[1]_\lambda^{2}[m]_\lambda^{2}/m$, and the weighted count of such curves is $[1]_\lambda^{2}[m]_\lambda^{2}/2$ after dividing by $\abs{\Aut\gamma}=2$ 
  
  \includegraphics{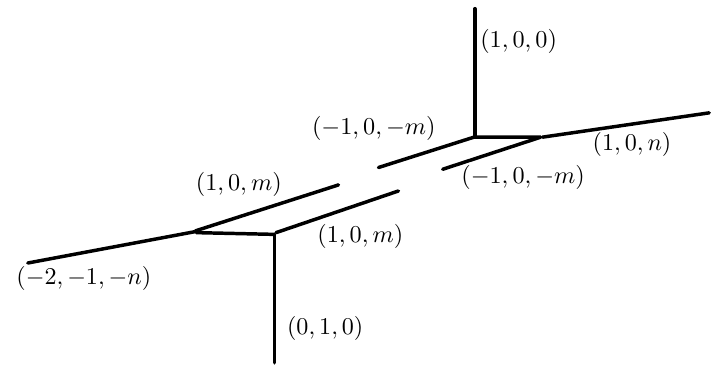}
  
  By equating our two different weighted counts of tropical curves, we therefore obtain the following formula.
  
  \begin{equation}\label{2m}2m[1]_\lambda^{2}=[2]_\lambda[2m]_\lambda/2+[1]_\lambda^{2}[m]_\lambda^{2}/2+ \sum_{a=1}^{m-1}[1]_\lambda^{2}[a]_\lambda^{2}\end{equation}
  
 Rearranging equations  (\ref{2m}) and (\ref{2m+1})  gives  the following two equations.
 \begin{equation}\label{odd}[2n+1]_\lambda=(2n+1)[1]_\lambda-[1]_\lambda\sum_{k=1}^{n}[k]_{\lambda}^{2}\end{equation}
\begin{equation}\label{even}[2n]_\lambda[2]_\lambda/2=[1]_\lambda^{2}(2n-[n]_\lambda^{2}/2-\sum_{k=1}^{n-1}[k]_{\lambda}^{2})\end{equation}

   Solving the above equations for $[n]_\lambda$ in terms of $[1]_\lambda$ gives 
 \begin{equation}\label{nq}[n]_\lambda=i^{-(n+1)}q^{n/2}+i^{n+1}q^{-n/2} \ \  \text{ where }\  q^{\pm\frac 12}:=\frac{-[1]_\lambda\pm\sqrt{[1]_\lambda^{2}-4}}2\ .\end{equation}
  
  \section{ Calculating absolute Gromov-Witten invariants of some toric manifolds}
  \label{calc section}
  Let $X$ be a smooth, compact,  $3$-complex dimensional toric manifold whose toric fan obeys the following convexity condition.

\begin{assumption}\label{ca} The intersection of any stratum of the fan of $X$ with the non-negative span of all other vectors in the fan of $X$ is $0$. \end{assumption}

This assumption has the consequence that any holomorphic curve contained in a boundary stratum of $X$ deforms into the interior of $X$. Without an assumption such as this, computation of Gromov-Witten invariants of $X$ necessarily involves interesting contributions from curves with components contained in the boundary of $X$. For an example of such contributions in the 2-dimensional setting, see \cite{tec}.

  For cohomology classes $\theta_{1},\dotsc,\theta_{k}\in H^{*}(X)$, and $\beta\in H_{2}(X)$, organize the corresponding Gromov-Witten invariants into a generating function as follows.
  \[\langle\theta_{1},\dotsc,\theta_{k} \rangle(X,\beta):=\sum_{g}\lambda^{2g-2+k}\int_{[\mathcal M_{g,k}(X,\beta)]]}\prod_{i=1}^{k}ev_{i}^{*}\theta_{i}\]
  In the above, $[\mathcal M_{g,k}(X,\beta)]$ indicates the virtual fundamental class of the moduli stack of connected genus $g$ curves in $X$ with $k$ labeled special points, representing  $\beta\in H_{2}X$. The tropical correspondence formula described below  holds when $\beta\neq 0$.
  
   For toric manifolds $X$ satisfying Assumption \ref{ca}, the intersection of $\beta$ with any complex codimension 1 toric  boundary stratum $S$ is a non-negative integer $d_{S}$ when $[\mathcal M(X,\beta)]$ is nonzero. Each such boundary stratum corresponds to an integral vector $\alpha_{S}$ in the fan of $X$. We need to count tropical curves with $d_{S}$ external edges  of derivative $\alpha_{S}$, and $k$ external edges of derivative $0$ (and no other external edges). Let $\bm\alpha$ be the corresponding list of derivatives of external edges. When $\beta\neq 0$, we have a correspondence formula in the following form:
  \begin{equation}\label{absolute correspondence}\langle\theta_{1},\dotsc,\theta_{k} \rangle(X,\beta)=W_{\bm\alpha}(\totb{\theta_{1}},\dotsc,\totb\theta_{k})\prod_{S}\frac{F^{d_{S}}}{d_{S}{!}}\end{equation}
  where $F=\lambda^{-1}+\sum_{g\geq 0} m_{g}\lambda^{2g-2+1}$ is some (relative) Gromov-Witten invariant,\footnote{The relative Gromow-Witten invariant $F$ is discussed in section \ref{ac section}. It can be regarded as a Gromov-Witten invariant of $\mathbb (CP^{1})^{3}$ relative to all but one boundary stratum.} and $W_{\alpha}(\totb\theta_{1},\dotsc,\totb\theta_{k})$ is a weighted count of tropical curves so that  the $i$th external edge with trivial derivative is constrained to $\totb{\theta_{i}}\subset{\mathbb R^{3}}$, constructed from $\theta_{i}$ as above equation (\ref{wdef}). For example, if all $\theta_{i}$ are  Poincare dual to a point, our weighted count of tropical curves counts tropical curves constrained to pass through $k$ points, and with $d_{S}$ external edges with derivative $\alpha_{S}$. To translate this into a weighted count $W_{\alpha}(\totb\theta)$ as in equation (\ref{wdef}), let $\theta$ be the product of the pullback of all $\theta_{i}$ to $\rh^{*}(\ex T^{3k+2\sum_{S} d_{S}})$. Then $W_{\bm\alpha}(\totb\theta_{1},\dotsc,\totb\theta_{k}):=W_{\bm\alpha}(\totb{\theta})$.   A more general form of this correspondence formula is proved in Section \ref{ac section}.  
   
   \

  For example, the toric fan of $(\mathbb CP^{1})^{3}$ satisfies Assumption \ref{ca}. Consider curves representing the homology class of the first $\mathbb CP^{1}$ constrained to pass through $1$ point. With the standard complex structure, there is a unique such genus $0$ curve, which contributes $\lambda^{-1}$ to our Gromov-Witten invariant. Theorem 3 from \cite{p1} implies that the higher genus contributions all vanish in this case, so our Gromov-Witten invariant is $\lambda^{-1}$. For our correspondence formula, we count tropical curves with $3$ external edges in directions $(\pm 1,0,0)$ and $0$ respectively, where the trivial external edge is constrained to a chosen point in $\mathbb R^{3}$. There is a unique such tropical curve, and it has weight $\lambda$. Our correspondence formula now reads
  \[\lambda^{-1}=\lambda F^{2}\]
  so $F=\lambda^{-1}$, and our correspondence formula simplifies to 
  \begin{equation}\label{absolute}\langle\theta_{1},\dotsc,\theta_{k} \rangle(X,\beta)=W_{\bm\alpha}(\totb{\theta_{1}},\dotsc,\totb\theta_{k})\prod_{S}(d_{S}!\lambda^{d_{S}})^{-1}\end{equation}
  
  A second example will pin down $[1]_\lambda$. Consider degree 1 curves in $\mathbb CP^{3}$ constrained to pass through $2$ points. With the standard complex structure, there is one such line passing through $2$ points, and in this case, Theorem 3 of \cite{p1} gives the following formula for our Gromov-Witten invariant.
  \[\lrb{\frac{\sin(\lambda/2)}{\lambda/2}}^{2}\] 
  For our correspondence formula, we want the weighted count of tropical curves with ends in the directions $(1,0,0)$, $(0,1,0)$, $(0,0,1)$, and $(-1,-1,-1)$, and two ends of derivative $0$ constrained to two generic points. There is one such tropical curve, which counts with weight $\lambda^{2}[1]_\lambda^{2}$. Our  correspondence formula then gives
  \[\lrb{\frac{\sin(\lambda/2)}{\lambda/2}}^{2}=\lambda^{2}[1]_\lambda^{2}\lambda^{-4}\]
  so
  \[\lrb{{2\sin(\lambda/2)}}^{2}=[1]_\lambda^{2}\ .\]
 Direct calculation in Section \ref{zero genus section} gives $[1]_\lambda=\lambda+\dotsb$, so we may remove the squares from the above equation, obtaining
 \begin{equation}[1]_\lambda=2\sin{\lambda/2}= i^{-2}(ie^{i\lambda/2})+i^{2}(ie^{i\lambda/2})^{-1}\end{equation}
 then equation (\ref{nq}) allows us to complete our correspondence formula.
 \begin{equation}[n]_\lambda=i^{-(n+1)}(ie^{i\lambda/2})^{n}+i^{n+1}(ie^{i\lambda/2})^{-n}=2\sin{n\lambda/2}\end{equation}

\subsection{Gromov-Witten invariants relative some toric boundary divisors}
\label{relative correspondence}

\

We can also write a tropical correspondence formula for Gromov-Witten invariants of a compact $3$-dimensional toric manifold $X$ relative to some chosen union of components of the toric boundary divisor of $X$. Call the vectors and cones in $X$ associated to this chosen divisor non-special. Call the vectors associated to toric boundary divisors not in our chosen divisor special, and any stratum of the toric fan of $X$ containing a special vector special. For our correspondence formula, we require the following combinatorial assumption.

\begin{assumption}\label{rca} The intersection of any special stratum of the fan of $ X$ with the non-negative span of the other vectors in the fan of $X$ is $0$. \end{assumption}

Let $\ex X$ be the explosion of $X$ relative  its chosen divisor. This exploded manifold $\ex X$ has a tropical part $\totb{\ex X}$  naturally identified with the subset of the fan of $\ex X$ that is non-special. Our correspondence formula below is for the Gromov-Witten invariants of $\ex X$, or equivalently,  the Gromov-Witten invariants of $X$ relative its chosen divisor.\footnote{ In the case that our divisor is smooth, we get the usual relative Gromov-Witten invariants, however we shall allow slightly more general cohomology classes to be used to define these invariants. }
Each curve in $\ex X$ has a tropical part that is a tropical curve in $\totb{\ex X}$. We can label the ends of this tropical curve and record their derivative as a list of vectors $\bm\alpha$. We can also record the degree $d_{S}$ of intersection with each special boundary divisor, labeled by a special vector $\alpha_{S}$. Let $[\mathcal M_{\bm d,\bm\alpha}(\ex X)]$ indicate the virtual fundamental class of the moduli stack of curves in $\ex X$ with this data. 

For each end of our curves with tropical part of derivative $\alpha$, there is a corresponding evaluation map $ev_{\alpha}$, with target an exploded manifold $\ex X_{\alpha}$ constructed as follows. If $\alpha=0$,  then $\ex X_{\alpha}:=\ex X$. If $\alpha\neq 0$, then $\alpha$ is an integral vector somewhere inside the image of a non-special cone in the fan of $ X$. Construct a $2$-dimensional fan by taking the projection with kernel spanned by $\alpha$ of all cones containing $\alpha$, and call the image of special vectors special. This fan is $\mathbb R^{k}$ times the fan of a $(2-k)$-dimensional toric manifold, with some divisor determined by non-special vectors. In this case, $\ex X_{\alpha}$ is $\ex T^{k}$ times the explosion of this toric manifold relative this divisor.  Let $\ex X_{\bm\alpha}$ be the product of $\ex X_{\alpha}$ for all vectors $\alpha$ in our list $\bm\alpha$ and let 
\[ev_{\bm\alpha}:[\mathcal M_{\bm d,\bm\alpha}(\ex X)]\longrightarrow \ex X_{\bm\alpha}\]
be the corresponding evaluation map.
 We are allowed to pull back and integrate cohomology classes from $\rh^{*}(\ex X_{\bm\alpha})$. This cohomology equals the inverse limit of the cohomology of toric manifolds with toric fans obtained from the fan corresponding to $\ex X_{\bm\alpha}$ by subdividing non-special strata. Any such cohomology class $\theta$ may be identified with a cohomology class in $\rh^{*}(\ex T^{N})$ to which we can associate a tropical constraint $\totb\theta$ as above equation (\ref{wdef}). Our correspondence formula is then
\begin{equation}\label{gac}\int_{[\mathcal M_{\bm d,\bm\alpha}]}\lambda^{2g-2+n}ev_{\bm\alpha}^{*}\theta=W_{\bm d,\bm\alpha}(\totb\theta)\prod_{S}\frac{F^{d_{S}}}{d_{S}!}\end{equation}
where $W_{\bm d,\bm\alpha}(\totb \theta)$ indicates a weighted count of tropical curves with $d_{S}$ ends with derivative $\alpha_{S}$ in addition to the ends corresponding to $\bm\alpha$, constrained by $\totb\theta$, and the product is over the special boundary strata $S$ of our toric manifold. As discussed in the previous section, $F$ is a Gromov-Witten invariant, and equal to $\lambda^{-1}$. This correspondence formula (\ref{gac}) subsumes both our previous tropical correspondence formulae, (\ref{correspondence}) and (\ref{absolute correspondence}).  We prove it in section \ref{ac section}.

\section{Donaldson-Thomas invariants}
 \label{DT section}
 
 Theorem 1 of \cite{gwdt} tells us that to calculate the reduced Donaldson-Thomas invariants of a $3$-dimensional toric manifold from the Gromov-Witten generating function counting possibly disconnected curves with no trivial connected components, we need to make the change of variables $q=-e^{i\lambda}$ and adjust by some factors of $\lambda$ and $q^{1/2}$. Our gluing formula for $F_{\gamma}$ implies that if $\gamma$ has $k$ ends with zero derivative, $F_{\gamma}/\lambda^{k}$ is some polynomial in  $q^{ 1/2}=ie^{i\lambda/2}$ and $q^{-1/2}$.  In particular, $F_{\gamma}/\lambda^{k}$ is an entire holomorphic function of $\lambda$. Any relation between these $F_{\gamma}$ obtained by calculating Gromov-Witten invariants in two different ways 
 leads to a polynomial relation in $q^{\pm 1/2}$, which is valid for all complex values of $\lambda$, and is therefore valid for  $q^{1/2}$ any complex number or formal variable. We may therefore define another invariant weighted count of tropical curves using $F^{DT}_{\gamma}$ in place of $F_{\gamma}$, where $F^{DT}_{\gamma}(q^{\frac 12})$ is a Laurent polynomial in $q^{\frac 12}$, and
 \[F^{DT}_{\gamma}(ie^{i\lambda/2})=F_{\gamma}/\lambda^{k}\]
  where $\gamma$ has $k$ edges with zero derivative. $F^{DT}_{\gamma}$ obeys the same gluing formula as $F_{\gamma}$, and is calculated the same way, except, when we see tropical curves $\gamma$ with $3$ external edges and no internal edges,   we use \[F_{\gamma}^{DT}=\frac{i^{-(n+1)}q^{n/2}+i^{(n+1)}q^{-n/2}}n\] when $\gamma$ has edges with derivative $\alpha$ and $\beta$ so that $\abs{\alpha\wedge\beta}=n$, and  we use $F_{\gamma}^{DT}=1$ when $\gamma$ has exactly one edge with zero derivative.
  
   Using $F^{DT}_{\gamma}$ in place of $F_{\gamma}$ in the formula (\ref{wdef}) for a weighted count of curves, but allowing possibly disconnected tropical curves with no trivial connected components give another weighted count of curves $W^{DT}_{\bm\alpha}$, which is some Laurent polynomial in $q^{\frac 12}$.\footnote{It would be nice to have a direct interpretation of $W_{\bm\alpha}^{DT}$ (times a normalizing factor depending on $\bm\alpha$) as a virtual count of appropriate ideal sheaves on $\ex T^{3}$. We will give no such interpretation in this paper. }

  Using the Gromov-Witten/Donaldson-Thomas correspondence from Theorem 1 in \cite{gwdt} and our correspondence formula (\ref{absolute}), we get the following correspondence formula for the reduced Donaldson-Thomas invariants of a compact toric manifold satisfying Assumption \ref{ca}.
 \[Z'_{DT}(X,q\vert \theta_{1},\dotsc,\theta_{k})_{\beta}=W^{DT}_{\bm\alpha}(\totb\theta_{1},\dotsc,\totb\theta_{k})\prod_{S}\frac{q^{d_{S}/2}}{d_{S}!}\]
 In the above, $Z'_{DT}$ indicates the reduced Donaldson-Thomas invariant using the notation of \cite{gwdt}. 
 
 \section{Interpretation using exploded manifolds}
\label{exploded section}

In the remainder of this paper, we explain how  the tropical correspondence and gluing formulae in the preceding sections all follow from the theory of Gromov-Witten invariants in exploded manifolds, and the tropical gluing formula for such Gromov-Witten invariants.

\subsection{$F _{\gamma}$}

\label{Fgamma section}

\

Let us define $F_{\gamma}$ as a kind of Gromov-Witten invariant, then verify that it satisfies the required correspondence formula, (\ref{correspondence}), and gluing formula, (\ref{gf}).  Cover the moduli stack of holomorphic curves in $\ex T^{3}$ by embedded Kuranishi charts $(\mathcal U,V,\hat f/G)$, as described in section 2.9 of \cite{evc}. Each such Kuranishi chart features a family  of (not-necessarily holomorphic) curves $\hat f$ with a finite group of automorphisms $G$, and  a $G$-equivariant, finite-dimensional obstruction bundle $V$ over the exploded manifold $\ex F(\hat f)$ parametrizing $\hat f$. Each of these families $\hat f$ has universal tropical structure, described in \cite{uts}. The tropical part, $\totb{\ex F(\hat f)}$  of $\ex F(\hat f)$ is largely determined by having universal tropical structure: in particular, Remark 3.3 of \cite{uts}, implies that for any curve $f$ in $\hat f$ with tropical part $\gamma$, the stratum of $\totb{\ex F(\hat f)}$ containing $f$ is isomorphic to $P_{\gamma}$, and the isomorphism is canonical once an isomorphism of the tropical part of $f$ with $\gamma$ is chosen. As we are in a special case where $\gamma$ is a curve in $\mathbb R^{3}$ and we are only considering automorphisms that fix the ends of $\gamma$, the automorphism group $G$ acts trivially on $P_{\gamma}$, and the isomorphism is canonical. With this identification, $\totb{ev}_{e}$ is the tropical part of the evaluation map $ev_{e}$ from $F(\hat f)$ which records the position of the external edge $e$.

The construction of a virtual fundamental class in \cite{vfc} using an embedded Kuranishi structure involves the intersection of weighted-branched perturbations of the section $\dbar:\ex F(\hat f)\longrightarrow V$ with $0$. Roughly speaking, this virtual fundamental class is locally a formal sum of rational weights times transverse intersections of sections of $V$ with $0$. Each such transverse intersection also has universal tropical structure, so the tropical part of any stratum containing a curve with tropical part $\gamma$ is canonically isomorphic to $P_{\gamma}$. The (real) dimension of these transverse intersections is twice the number of external edges of $\gamma$, so the maximal dimension of their tropical part is the number of external edges of $\gamma$. It follows that $\gamma$ is general in the sense of Definition \ref{general def} if and only if the corresponding  strata have maximal tropical dimension, and have maximal dimensional (tropical) image under the evaluation map that records the position of  external edges.\footnote{Each stratum of an exploded manifold projects down to an underlying manifold, called its smooth part, whose real dimension is the real dimension of the exploded manifold minus twice the tropical dimension of the stratum. In our case, a stratum of the moduli stack with maximal tropical dimension has $0$-dimensional smooth part. The smooth part of the moduli stack of curves likely corresponds  to the reader's natural conception of what the moduli stack of curves should be, and from this perspective, a stratum with maximal tropical dimension looks like a maximally degenerate stratum. In this paper, we are making all our calculations in this apparently degenerate setting. A proof of our results using standard geometry would involve systematically keeping track of many blowups.}

 As explained in section 7 of \cite{vfc}, the contribution of curves with tropical part $\gamma$ to the virtual fundamental class is formalized using the tropical completion $[\mathcal M\tc \gamma]$ of the virtual fundamental class at $\gamma$. If $\gamma$ is general, the stratum of $[\mathcal M_{\alpha}]$ containing curves with tropical part $\gamma$ has maximal tropical dimension (equal to the complex virtual dimension of the moduli stack), so  $[\mathcal M\tc \gamma]$ is some formal weighted sum of exploded manifolds isomorphic to $\ex T^{n}$.  The tropical completion $ev_{\bm\alpha}\tc \gamma$ of the evaluation map $ev_{\bm\alpha}$ restricted to this strata is some map $ev_{\bm\alpha}\tc\gamma:\ex T^{n}\longrightarrow \ex T^{2n+k}$ whose tropical part is the extension  of $\totb{ev}_{\bm\alpha}$ to an affine  map. All such maps  $ev_{\bm\alpha}\tc\gamma$ differ only by the translations $\ex T^{2n+k}\longrightarrow\ex T^{2n+k}$ induced by multiplying coordinates by constants. We can therefore choose a differential form $\theta_{\bm\gamma}\in\rh^{2n}(\ex T^{2n+k})$ so that the integral of the pullback of $\theta_{\bm\gamma}$ under any such map is $1$. Define $F_{\gamma}$ as the following Gromov-Witten invariant,
 \begin{equation}\label{Fdef}F_{\gamma}:=\abs{\Aut \gamma}\int_{[\mathcal M_{\bm\alpha}\tc\gamma]}\lambda ^{2g-2+n}(ev_{\bm\alpha}\tc\gamma)^{*}\theta_{\gamma}\end{equation}
 where $n$ is the number of $3$-vectors in $\bm\alpha$ (or the number of ends of the curves we are studying) and $g$ is the function on $[\mathcal M_{\bm\alpha}]$ recording the genus of curves. The extra factor of $\abs{\Aut \gamma}$ appears because we want to apply a gluing formula which is valid for the  moduli stack of curves with a chosen isomorphism of their tropical part with $\gamma$.
 
 \subsection{Correspondence formula for  Gromov-Witten invariants of $\ex T^{3}$}\label{correspondence section}
 
 \
 
  We now verify that our correspondence formula (\ref{correspondence}) holds with $F_{\gamma}$ defined as above.  Lemma 7.7 of \cite{vfc} gives that for any closed  $\theta\in\rof^{2n} (\ex T^{2n+k})$, 
 \begin{equation}\label{tc}\int_{[\mathcal M_{\bm\alpha}]}ev_{\bm\alpha}^{*}\theta=\sum_{\gamma}\int_{[\mathcal M_{\bm\alpha}\tc \gamma]} (ev_{\bm\alpha}^{*}\theta)\tc\gamma\ .\end{equation}
In particular, let $\theta'$ be a closed differential form on $X$, a compact $(2n+k)$-dimensional toric manifold. Define $\theta\in\rof^{2n}(\ex T^{2n+k})$ by  pulling back $\theta'$ to $\expl X$, then pushing it forward to $\ex T^{2n+k}$ using a refinement map $\expl X\longrightarrow \ex T^{2n+k}$. Equation (\ref{tc}) above is valid for any such $\theta$. Let $\totb\theta$ be as discussed above equation (\ref{wdef}). Even with $\totb \theta$ shifted generically, we can choose our refinement map $\expl X\longrightarrow \ex T^{2n+k}$ so that the strata of $\totb\theta$ are strata of the corresponding subdivision of $\totb{\ex T^{2n+k}}=\mathbb R^{2n+k}$. If the image of $\gamma$ is not in $\totb\theta$, then $\int(ev_{\bm\alpha}^{*}\theta)\tc\gamma=0$, so as with equation (\ref{wdef}), the right hand side of equation (\ref{tc}) is a sum over $\gamma$ so that $\totb{ev}_{\bm\alpha}(\gamma)$ is in $\totb\theta$. Using the notation of equation (\ref{wdef}), when $\totb{ev}_{\bm\alpha}\gamma\in\totb\theta_{S}$, then $\theta\tc{(\totb{ev}_{\alpha}\gamma)}$ is equal to $m_{S}$ times the Poincare dual to the $\ex T^{n+k}\subset \ex T^{2n+k}$ with tropical part the extension of $\totb\theta_{S}$ to an affine subspace. Moreover, the integral of the pullback of this using any of the maps $\ex T^{n}\longrightarrow \ex T^{2n+k}$ appearing in $[\Msw_{\bm\alpha}\tc\gamma]$ is the $\abs{\mathbb Z^{2n+k}/ev_{\alpha}(\Lambda_{\gamma})\oplus\Lambda_{\totb\theta_{S}}}$ appearing in equation (\ref{wdef}). Using that $(ev_{\bm\alpha}^{*}\theta)\tc\gamma=(ev_{\bm\alpha}\tc\gamma)^{*}\theta\tc{(\totb{ev}_{\bm\alpha}\gamma)}$ together with our definition of $F_{\gamma}$, and equation (\ref{tc}), we obtain
\[\int_{\mathcal M_{\bm\alpha}}\lambda^{2g-2+n}ev^{*}_{\bm\alpha}\theta=\sum_{S}\sum_{\gamma\in\totb{ev}_{\bm\alpha}^{-1}\totb\theta_{S}}\frac{F_{\gamma}}{\Aut\gamma}m_{S}\abs{\mathbb Z^{2n+k}/ev_{\alpha}(\Lambda_{\gamma})\oplus\Lambda_{\totb\theta_{S}}}\]
which is our correspondence formula (\ref{correspondence}) using equation (\ref{wdef}) to define $W_{\bm\alpha}(\totb\theta)$. The correspondence formula using equation (\ref{W}) to define $W_{\bm\alpha}(\totb\theta)$ is a special case.

\subsection{Gluing formula}
\label{gf section}

\

We shall now verify that $F_{\gamma}$ satisfies the claimed gluing formula, (\ref{gf}). The gluing formula, Theorem 4.7 of \cite{egw},\footnote{The gluing formula from \cite{egw} uses an old construction of the virtual fundamental class. Its proof is similar, but a bit easier using the more recent construction of virtual fundamental class from  \cite{evc,vfc}. I will publish an updated, more general version of this gluing formula shortly.} gives a formula for $F_{\gamma}$ as follows. Let $\bm\alpha_{v}$ be the list of the derivatives of edges of $\gamma$ leaving $v\in V(\gamma)$. As $\gamma$ is not general in our sense whenever it has internal edges with $0$ derivative, we may assume that all internal edges of $\gamma$ have nonzero derivative. (The gluing formula from \cite{egw} also implies that such tropical curves do not contribute to our Gromov-Witten invariants.) Suppose that $\gamma$ has $n'$ internal edges and $n$ external edges, $k$ of which have $0$ derivative. We shall express $F_{\gamma}$ as an integral over $\ex T^{2n+k}\times \ex T^{2n'}$, where the extra factors of $\ex T^{2}$  are  the quotient of $\ex T^{3}$ by the action of weight $\alpha_{e}/\abs{\alpha_{e}}$ for internal edges $e$.  Let $\eta_{v}$ be the pullback to $\ex T^{2n+k+2n'}$ of the pushforward\footnote{See section 5.3 of \cite{vfc} for the definition of pushing forward a differential form using an evaluation map such as $ev_{\bm\alpha}$.} of $\lambda^{2g-2+n_{v}}$ using $ev_{\bm\alpha_{v}}$, and let $\theta_{\gamma}$ be the pullback to $\ex T^{2n+k+2n'}$ of the $\theta_{\gamma}$ used to define $F_{\gamma}$ in equation (\ref{Fdef}). Then Theorem 4.7 of \cite{egw} implies
\begin{equation}\label{kgamma}F_{\gamma}=k_{\gamma}\int_{\ex T^{2n+k+2n'}}\theta_{\gamma}\bigwedge_{v}\eta_{v}\text{ where }k_{\gamma}:=\frac1{\abs{\Aut\gamma}}\prod_{e\in \ie(\gamma)}\abs{\alpha_{e}}\ .\end{equation}

 Let us use our correspondence formula on the above expression. Our correspondence formula also works (with the same proof) for disconnected curves. In particular, consider curves with one connected component with data $\bm\alpha_{v}$ for each $v\in V(\gamma)$. See section 6 of \cite{vfc} for a discussion of the technical details involved in proving that the virtual fundamental class of the moduli stack of such disconnected curves is $\prod_{v}[\Mod_{\bm{\alpha_{v}}}]$. For a tropical curve $\coprod_{v}\gamma_{v}$ with connected components $\gamma_{v}$, this  implies that $F_{\coprod_{v}\gamma_{v}}=\prod_{v}F_{\gamma_{v}}$. Use 
\[ev_{\{\bm\alpha_{v}\}}:\prod_{v\in V(\gamma)}[\mathcal M_{\bm\alpha_{v}}]\longrightarrow \ex T^{2n+k+4n'}\]
to indicate the product of  $ev_{\bm\alpha_{v}}$ for all $v$. Consider $\ex T^{2n+k+4n'}$ as the target, $\ex T^{2n+k}$, of $ev_{\alpha}$, times two copies of the evaluation space $\ex T^{2}$ for each internal edge of $\gamma$, as each such edge corresponds to two edges in $\{\bm\alpha_{v}\}$.   Let $\theta_{\Delta}$ be the Poincare dual of the diagonal inclusion  $\ex T^{2n+k+2n'}\subset \ex T^{2n+k+4n'}$, and use $\theta_{\gamma}$ to indicate the pullback of $\theta_{\gamma}$ from $\ex T^{2n+k}$ to $\ex T^{2n+k+4n'}$. Then  equation (\ref{kgamma}) above is equivalent to the following expression.
\[F_{\gamma}=k_{\gamma}\int_{\prod_{v}[\mathcal M_{\bm\alpha_{v}}]}ev_{\{\bm\alpha_{v}\}}^{*}(\theta_{\gamma}\wedge\theta_{\Delta})\]
Applying our correspondence formula to this Gromov-Witten invariant gives
\begin{equation}\label{c4}F_{\gamma}=k_{\gamma}W_{\{\bm\alpha_{v}\}}(\totb\theta_{\gamma}\cap \totb\theta_{\Delta}) \end{equation}
where $W_{\{\bm\alpha_{v}\}}$ indicates a weighted count of tropical curves with connected components having data $\bm\alpha_{v}$, and the constraint $\totb\theta_{\gamma}\cap\totb\theta_{\Delta}$ means the construction above equation (\ref{wdef}) applied to $\theta_{\gamma}\wedge\theta_{\Delta}$.  We can think of this as a constraint $\totb\theta_{\gamma}\subset \mathbb R^{2n+k}$ on the edges corresponding to the external edges of $\gamma$ and a separate constraint $\totb\theta_{\Delta}\subset\mathbb R^{4n'}$ on those edges corresponding to  interior edges of $\gamma$. For our correspondence formula to work, we must choose this constraint  $\totb\theta_{\Delta}$ some generic shift of the diagonal. The constraint $\totb\theta_{\gamma}$ can be regarded as the constraint corresponding to a linear subspace  $\totb\theta_{0}$ divided by the index $m_{0}$ of $\totb{ev}_{\bm\alpha}(\Lambda_{\gamma})\oplus\Lambda_{\totb\theta_{0}}$,  so $m_{0}\totb\theta_{\gamma}=\totb\theta_{0}$ and our correspondence formula tells us 
\[m_{0}=\abs{\mathbb Z^{2n+k}/\totb{ev}_{\bm\alpha}(\Lambda_{\gamma})\oplus\Lambda_{\totb\theta_{0}}}\ .\]
 Then equation (\ref{wdef}) gives the following expression for $W_{\{\bm\alpha_{v}\}}(\totb\theta_{\gamma}\cap\totb\theta_{\Delta})$:
\begin{equation}\label{c1}\sum_{\gamma'\subset\totb{ev}_{\{\bm\alpha_{v}\}}^{-1}(\totb\theta_{0}\cap\totb\theta_{\Delta})} \frac1{m_{0}}\abs{\mathbb Z^{2n+k+4n'}/\oplus_{v}\totb{ev}_{\bm\alpha_{v}}(\Lambda_{\gamma'_{v}})\oplus\Lambda_{\totb\theta_{0}}\oplus \Lambda_{\totb\theta_{\Delta}}}\prod_{v}\frac{F_{\gamma'_{v}}}{\abs{\Aut\gamma'_{v}}} \end{equation}
where the sum is over disconnected tropical curves $\gamma'$ with components $\gamma'_{v}$ satisfying our constraints on edges. We can break up the size of the quotient of $\mathbb Z^{2n+k+4n'}$ above into two factors as follows.
\begin{equation}\label{c2}\abs{\mathbb Z^{2n'}/h(\oplus_{v}\Lambda_{\gamma'_{v}})}\abs{\mathbb Z^{2n+k}/g(\ker h)\oplus \Lambda_{\totb \theta_{0}}}\end{equation}
where $h$ is the composition of $\prod_{v}\totb{ev}_{\bm\alpha_{v}}$ with the projection to $\mathbb Z^{2n'}$ with kernel $\mathbb Z^{2n+k}\oplus \Lambda_{\totb\theta_{\Delta}}$, and $g$ is the composition of $\prod_{v}\totb{ev}_{\bm\alpha_{v}}$ with projection to $\mathbb Z^{2n+k}$. We must compare these to the factor $m_{\gamma'}:=\abs{\mathbb Z^{3n'}/A_{\gamma'}(\oplus_{v}\Lambda_{\gamma_{v}}\oplus\mathbb Z^{n'})}$ appearing in equation (\ref{gf}). Applying the snake lemma to the following commutative diagram of short exact sequences
\[\begin{tikzcd}\mathbb Z^{n'}\dar\rar&\oplus_{v}\Lambda_{\gamma_{v}}\oplus \mathbb Z^{n}\rar\dar{A_{\gamma'}} &\oplus_{v}\Lambda_{v}\dar{h}
\\ \mathbb Z^{n'}\rar&\mathbb Z^{3n'}\rar&\mathbb Z^{2n'}\end{tikzcd}\]
gives the following long exact sequence.
\[\begin{tikzcd}0\rar&\ker A_{\gamma'}\rar&\ker h\rar&\prod_{e\in\ie(\gamma)}\mathbb Z_{\abs{\alpha_{e}}}\rar&\mathbb Z^{3n'}/ A_{\gamma'}\rar&\mathbb Z^{2n'}/h\rar&0\end{tikzcd}\]
Therefore,
\begin{equation}\label{c3}\abs{\ker h/\ker A_{\gamma'}}m_{\gamma'}=\abs{\mathbb Z^{2n'}/h(\oplus_{v}\Lambda_{\gamma_{v}})}\prod_{e\in\ie(\gamma)}\abs{\alpha_{e}}\ .\end{equation}
Inspection of the equations (\ref{equation}) and (\ref{delta equation}) used to define $P_{\gamma}$ and $A_{\gamma'}$ respectively reveals a natural free action of $P_{\gamma}$  on the kernel of $A_{\gamma}$. The fact that $\gamma$ is general implies that $P_{\gamma}$ has the same dimension as $A_{\gamma}$ and  that $\Lambda_{\gamma}=\ker A_{\gamma'}$. Both $\totb{ev}_{\bm\alpha}$ and $g$ are injective, therefore 
\begin{equation}\label{c0}\abs{\ker h/\ker A_{\gamma'}}\abs{\mathbb Z^{2n+k}/g(\ker h)\oplus\Lambda_{\totb\theta_{0}}}=\abs{\mathbb Z^{2n+k}/\totb{ev}_{\bm\alpha}(\Lambda_{\gamma})\oplus\Lambda_{\totb\theta_{0}}}=m_{0}\ .\end{equation}
Note that solutions $\gamma'$ to $\totb{ev}_{\{\bm\alpha_{v}\}}(\gamma')\in \totb\theta_{0}\cap \totb\theta_{\Delta}$ correspond to $\delta$-deformed tropical curves $\gamma'$ with edges constrained by $\totb\theta_{0}$, where $\delta$ depends on our choice of translation of $\totb{\theta}_{\Delta}$. This observation together with  equations (\ref{c4}), (\ref{c1}), (\ref{c2}), (\ref{c3}), and (\ref{c0}) together imply our gluing formula
\[F_{\gamma}=\sum_{\gamma'}m_{\gamma'}\prod_{v}\frac{F_{\gamma_{v}'}}{\abs{\Aut\gamma_{v}}}\]
which is equation (\ref{gf}).
\subsection{Correspondence formula for absolute invariants and invariants relative to some toric boundary divisors.}\label{ac section}

\

 Let $X$ be a compact $3$-dimensional toric manifold, with a divisor $D$ consisting of some collection of irreducible components of the boundary divisor of $X$. Let $\ex X:=\expl (X,D)$. Our correspondence formula (\ref{gac}) gives Gromov-Witten invariants of $\ex X$ in terms of the weighted counts of tropical curves corresponding to Gromov-Witten invariants of $\ex T^{3}$. Let $\ex Y$ be the explosion of $X$ relative  its toric boundary divisor. $\ex Y$ is a refinement of $\ex T^{3}$, and its Gromov-Witten invariants are naturally identified with the Gromov-Witten invariants of $\ex T^{3}$.\footnote{For a proof that refinement does not affect Gromov-Witten invariants of exploded manifolds, see section 5.2 of \cite{egw}. In the log setting, see \cite{ilgw}.} Exploding the identity map from $X$ with its toric boundary divisor to $X$ with $D$ gives a map $\ex Y\longrightarrow \ex X$. As we already know that the Gromov-Witten invariants of $\ex Y$ obey our tropical correspondence formula,   let us reformulate equation (\ref{gac}) as giving the Gromov-Witten invariants of $\ex X$ in terms of the Gromov-Witten invariants of $\ex Y$. Let $(\bm d,\bm \alpha)$ indicate a list including the vectors in $\bm \alpha$, and $d_{S}$ copies of the vector $\alpha_{S}$ for special vectors $\alpha_{S}$ in the fan of $X$ indicating toric boundary components not contained in $D$. Our Gromov-Witten invariants of $\ex Y$ use an evaluation map
 \[ev_{(\bm d,\bm\alpha)}:[\Msw_{(\bm d,\bm\alpha)}(\ex Y)]\longrightarrow \ex Y_{(\bm d,\bm\alpha)}\]
that keeps track of the location of all the ends of curves. The corresponding evaluation map for the Gromov-Witten invariants of $\ex X$
\[ev_{\bm\alpha}:[\mathcal M_{\bm d,\bm\alpha}(\ex X)]\longrightarrow \ex X_{\bm\alpha}\]
only keeps track of ends labeled by vectors in $\bm\alpha$, because the vectors in direction $\alpha_{S}$ no longer correspond to ends of curves in $\ex X$. The map $\ex Y\longrightarrow \ex X$ induces a natural map 
\[\pi:\ex Y_{(\bm d,\bm\alpha)}\longrightarrow \ex X_{\bm\alpha}\ .\]
With this, we can reformulate equation (\ref{gac}) as a statement that the Gromov-Witten invariants of $\ex Y$ are given by our weighted counts of tropical curves, and that the Gromov-Witten invariants of $\ex X$ are determined by the Gromov-Witten invariants of $\ex Y$ as follows:
\begin{equation}\label{eta c}(ev_{\bm\alpha})_{!}\lambda^{2g-2+n}=\pi_{!}(ev_{(\bm d,\bm\alpha)})_{!}\lambda^{2g-2+n}\prod_{S}\frac{F^{d_{S}}}{d_{S}!}\end{equation}
where $(ev_{\bm\alpha})_{!}\lambda^{2g-2+n}$ is a laurent series in $\lambda$ with coefficients in $\rh^{*}(\ex X_{\bm\alpha}) $ uniquely determined by the property that 
\[\int _{[\mathcal M_{\bm d,\bm\alpha}(\ex X)]}\lambda^{2g-2+n}ev_{\bm\alpha}^{*}\theta=\int_{\ex X_{\bm\alpha}}\theta\wedge(ev_{\bm\alpha})_{!}\lambda^{2g-2+n}\]
for all $\theta\in \rh^{*}(\ex X_{\bm\alpha})$. As with equation (\ref{gac}), equation (\ref{eta c}) only applies under assumptions on the divisor $D$ and the toric fan of $X$ such as Assumption \ref{rca}.

\

Before discussing the general case, consider the following special case of our correspondence formula (\ref{gac}).  Let $\ex X$ be $\ex T^{2}\times\expl(\mathbb CP^{1},0)$.\footnote{Although such an $\ex X$ is not exactly in the form covered by Assumption \ref{rca}, it has a refinement that is in the form covered by Assumption \ref{rca} so we can expect our correspondence formula (\ref{gac}) to hold.} The tropical part of $\ex X$, $\totb{\ex X}$ is $\mathbb R^{2}\times [0,\infty)$. Any holomorphic curve in $\ex X$ has tropical part a tropical curve in $\totb{\ex X}$ which obeys the usual balancing condition in the interior, and has the sum of all derivatives leaving a vertex on the boundary some non-negative multiple of $(0,0,1)$. 

 Consider the case of a tropical curve with a single vertex on the boundary and  a single edge, of derivative $(0,0,1)$.  The evaluation map at this edge has target the $4$-dimensional exploded manifold $\ex T^{2}$. The degree of this evaluation map restricted to the curves of genus $g$ is a Gromov-Witten invariant, $n_{g}$. Make the definition
\[F:=\sum_{g=0}^{\infty}n_{g}\lambda^{2g-2+1}\ .\]
It is elementary to check that $n_{0}=1$, so $F=\lambda^{-1}+\dotsc$. We shall prove that the above is the $F$ appearing in formula (\ref{eta c}), (\ref{gac}) and the special case (\ref{absolute correspondence}), from which it follows that $F=\lambda^{-1}$. 

In this special case, where $\bm{\alpha}=\{(0,0,1)\}$, our correspondence formula (\ref{gac}) holds by the  definition of $F$. We shall use this fact to prove that our correspondence formula holds in general.

\

To relate our Gromov-Witten invariants of $\ex T^{3}$ to the Gromov-Witten invariants of a $3$-dimensional toric manifold $X$ relative to some (possibly empty) divisor $D$ consisting of some number of irreducible components of the toric boundary divisor, we shall choose an appropriate degeneration of $X$, or more accurately a degeneration of the explosion of $(X,D)$.  Our degeneration will correspond to a subdivision of the moment polytope of $X$ dual to the barycentric subdivision.  Let $\hat X$ be a $4$-dimensional toric manifold with fan containing cones spanned by the following sets of vectors for each cone in the fan of $ X$ spanned by $\beta_{1}$, $\beta_{2}$ and $\beta_{3}$:
\begin{itemize}
\item $\{(0,\beta_{1}),(0,\beta_{2}),(0,\beta_{3}),(1,\sum\beta_{i})\}$
\item $\{((0,\beta_{1}),(0,\beta_{2}),(1,\beta_{1}+\beta_{2}),(1,\sum\beta_{i}))\}$
\item $\{(0,\beta_{1}),(1,\beta_{1}),(1,\beta_{1}+\beta_{2}),(1,\sum\beta_{i})\}$
\item $\{(1,0),(1,\beta_{1}),(1,\beta_{1}+\beta_{2}),(1,\sum\beta_{i})\}$
\end{itemize}
and all permutations of the roles of $\beta_{i}$.  Projecting these fans onto the first coordinate corresponds to a toric map  $\pi:\hat X\longrightarrow \mathbb C$. This $\pi$ is a normal crossing degeneration with generic fiber $X$. Use the adjective `special' to describe any vector in the toric fan of $X$ corresponding to an irreducible component of the boundary divisor not in $D$. The special vectors $\beta$ in the fan of $X$ give us special vectors, $(0,\beta)$ in the fan of $\hat X$ and a divisor $\hat D$ consisting of all non-special components of the toric boundary divisor. Consider  $\mathbb C$ to have the divisor $0$. Then, as explained in \cite{iec}, we can apply the explosion functor to $(\hat X,\hat D)$ to obtain an exploded manifold $\hat {\ex X}$ which is the total space of a smooth family of exploded manifolds.
\[\expl \pi:\hat{\ex X}\longrightarrow \expl (\mathbb C,0)=\et 1{[0,\infty)}\]
As proved in \cite{vfc}, Gromov-Witten invariants do not change in such families of exploded manifolds, so we can compute the Gromov-Witten invariants of the generic fiber, $\expl( X,D)$, using the Gromov-Witten invariants of a more interesting fiber $\ex X:=(\expl \pi)^{-1}(1\e{1})$ (using notation from \cite{iec}). The tropical part $\totb{\hat {\ex X}}$ of $\hat{\ex X}$ is naturally identified with the subset of the fan of $\hat X$ spanned by the non-special vectors above, and the tropical part $\totb {\ex X}$ of $\ex X$ is the subset of this which is the inverse image of $1$. 

Let us apply the gluing formula, Theorem 4.7 of \cite{egw}, to $\ex X$.
This formula involves Gromov-Witten invariants $\eta_{v}$ of exploded manifolds $\ex X\tc v$ obtained from $\ex X$ using tropical completion (described in section 7 of \cite{vfc} as well as in \cite{egw}.) These exploded manifolds $\ex X\tc v$ can be constructed as follows:
Consider $(1,v)$ as a point somewhere inside the (non-special part of the) toric fan of $\hat X$. For each cone $\sigma$ of the toric fan of $\hat X$ that contains $(1,v)$, consider the projection, $\sigma_{v}\subset\mathbb R^{3}$ of $\sigma$ using the projection with kernel containing $(1,v)$ and with image the plane with first coordinate $0$. The union of such $\sigma_{v}$ determines a subdivison of $\mathbb R^{3}$. The image of all special vectors in such a $\sigma$ determine special vectors for this subdivision. If $v$ is a $0$-dimensional stratum of $\totb{\ex X}$, this subdivision is the fan of a toric manifold with divisor specified by the non-special strata. $\ex X\tc v$ is the explosion of this toric manifold relative to this divisor. The Gromov-Witten invariants of $\ex X\tc v$ can be regarded as the Gromov-Witten invariants of this toric manifold relative this divisor. If $v$ is in a $d$-dimensional stratum, this subdivision is $\mathbb R^{d}$ times a toric fan for a $(3-d)$-dimensional toric manifold with a divisor determined by the non-special strata. In this case, $\ex X\tc v$ is $\ex T^{d}$ times the explosion of this toric manifold  relative to this divisor.

Let us examine $\ex X\tc v$ more closely in the case that $v$ is a $0$-dimensional stratum of $\totb{\ex X}$. Such a $v$ is  in the form $\sum_{i=1}^{k}\beta_{i}$ for $k=0,1,2,$ or $3$. Some number $d\leq k$ of these $\beta_{i}$ will be special; without losing generality, suppose that $\beta_{1},\dotsc,\beta_{d}$ are special. If $d=0$,   $\ex X\tc v$ is a refinement of $\ex T^{3}$, and our correspondence formula  (together with the fact, proved in section 5.2 of \cite{egw}, that refinement does not affect Gromov-Witten invariants)  tells us the Gromov-Witten invariants of $\ex X\tc v$. Otherwise, the non-special strata are some subdivision of the non-negative span of $-\beta_{1},\dotsc,-\beta_{k}$, $\beta_{d+1},\dotsc,\beta_{k}$ and  all vectors that span a cone in the fan of  $X$ along with $\beta_{1},\dotsc,\beta_{k}$.  Our toric manifold minus its divisor (determined by non-special vectors) is $\mathbb C^{d}\times(\mathbb C^{*})^{3-d}$. The virtual dimension of the moduli space of holomorphic curves in such an $\ex X\tc v$ is readily computed as follows. Each holomorphic curve in $\ex X\tc v$ has a tropical part, consisting of a tropical curve in the non-special strata of our fan. Let $\alpha_{1},\dotsc,\alpha_{n}$ be the outgoing derivatives of the ends of such a tropical curve. Specifying $\alpha_{i}$ determines the homology class represented by our curve, and therefore the virtual dimension of the corresponding moduli space of holomorphic curves.  There exist unique integers $n_{1},\dotsc,n_{d}$ so that the following balancing condition holds. 
\[n_{1}\beta_{1}+\dotsb+n_{d}\beta_{d}+\sum_{i=1}^{n}\alpha_{i}=0\]
Then the (real) virtual dimension of the corresponding moduli space of curves is 
\begin{equation}\label{vdim}2(\sum_{i=1}^{d}n_{i}+n)\ .\end{equation}
This is the point where Assumption \ref{rca} on the toric fan of $ X$ comes in useful.\footnote{When Assumption \ref{ca} does not hold, the Gromov-Witten invariants of $\ex X\tc v$ become more interesting, and enter into the correspondence formula in an essential way.}  Each $\alpha_{i}$ must be in the non-negative span of $-\beta_{1},\dotsc,-\beta_{k},\beta_{d+1},\dotsc,\beta_{k}$ and vectors that span a cone in the fan of  $ X$ along with $\beta_{1},\dotsc,\beta_{k}$. Assumption \ref{rca}  implies that each $n_{i}\geq 0$. Therefore, the above virtual dimension (\ref{vdim}) is the virtual dimension of the space of curves in $\ex T^{3}$ corresponding to a tropical curve created by adding external edges with derivative $\beta_{1},\dotsc,\beta_{d}$ to our tropical curve in $\totb{\ex X}\tc v$ until it is balanced. 

The same formula holds for the virtual dimension of the moduli stack of curves in any of our $\ex X\tc v$. Our gluing formula, Theorem 4.7 of \cite{egw}, then implies that there is a similar formula for the virtual dimension of the moduli stack of curves containing a curve with tropical part $\gamma$ in $\totb{\ex X}$: there is a unique way to complete this tropical curve to a balanced tropical curve in (a subdivision of) $\mathbb R^{3}$ by adding ends with special derivative; the  (real) virtual dimension of the corresponding moduli stack of curves is equal to twice the number of ends of this balanced tropical curve. Define the tropical virtual dimension to be the number of ends of this balanced tropical curve. Define a tropical evaluation map $\totb{ev}_{\bm\alpha}$ similarly to the case of tropical curves in $\mathbb R^{3}$, and similarly to  Definition \ref{general def}, call a tropical curve general if its  space of deformations has dimension equal to the tropical virtual dimension and $\totb {ev}_{\bm\alpha}$ is injective on this space. As in the case of $\ex T^{3}$, we shall  show that we  need only consider the contribution of general tropical curves, and we shall compute such contributions using Theorem 4.7 of \cite{egw}. 

Say that our correspondence formula $(\ref{gac})$ holds at a vertex $v$ of a tropical curve  if it holds for $\ex X\tc v$ and $\bm\alpha$ the list of derivatives of edges leaving $v$.

\begin{claim}\label{induct} Suppose that for every general tropical curve in $\totb {\ex X}$ with  exterior edge derivatives $\bm{\alpha}$, our correspondence formula holds for all vertices on the boundary of $\totb{\ex X}$. Then our correspondence formula also holds for $\bm{\alpha}$.
\end{claim}

\

To prove Claim \ref{induct}, it shall be helpful for us to use $\ex Y$ to indicate the exploded manifold constructed using the same construction as $\ex X$, but using the entire toric boundary divisor instead of just $D$. As noted earlier, there is a natural map $\ex Y\longrightarrow \ex X$. $\ex Y$ is some refinement of $\ex T^{3}$ determined by a subdivision $\totb{\ex Y}$ of $\mathbb R^{3}$. The tropical part, $\totb{\ex X}$ of $\ex X$,  can be regarded as some union of `non-special' strata of $\totb {\ex Y}$.

We need to choose a representative of our our cohomology class  on our evaluation space $\ex X_{\bm\alpha}$. Recall, from above equation (\ref{gac}), that $\ex X_{\bm\alpha}=\prod_{e\in \ee(\gamma)}\ex X_{\alpha_{e}}$. As we are in a slightly more general situation to that described above equation (\ref{gac}), let us describe $\ex X_{\alpha_{e}}$ in detail. If $\alpha_{e}=0$, then $\ex X_{\alpha_{e}}=\ex X$. Otherwise,  consider the $3$-dimensional toric manifold $\hat X_{\alpha_{e}}$ with fan the projection, with kernel spanned by $\alpha_{e}$,  of all cones in the toric fan of $\hat X$ containing $\alpha_{e}$, and with special vectors the image of special vectors in these cones under this projection. Like $\hat X$, $\hat X_{\alpha_{e}}$ has a divisor that is the toric boundary divisor minus special components, and also comes with a toric map to $\mathbb C$.  Exploding this map gives a family of exploded manifolds, and $\ex X_{\alpha_{e}}$ is the fiber over $1\e 1$.  (As in the cases discussed earlier, when this projected subdivision is $\mathbb R^{d}$ times the fan of a toric manifold, we use $\ex T^{d}$ times the explosion of this toric manifold.)

We can make the same construction for $\ex Y$. There is a natural map $\ex Y_{\alpha_{e}}\longrightarrow \ex X_{\alpha_{e}}$, however we want a tropical  curve $\gamma$ in $\totb{\ex X}$ to correspond to the unique tropical curve $\gamma'$ in $\totb{\ex Y}$ obtained by adding edges with derivatives special vectors until the resulting curve is balanced. If $(\bm d,\bm\alpha)$ indicates this list of derivatives of external edges of $\gamma'$, there is a natural map
\[\pi:\ex Y_{(\bm d,\bm\alpha)}\longrightarrow \ex X_{\bm\alpha}\]
which forgets the extra edges in $\gamma'$. We can choose a representative $\theta$ for our cohomology class in $\rhf^{2m}(\ex X_{\bm\alpha})$ with support $\totb \theta$ in $\totb{\ex X_{\bm\alpha}}$ having codimension $m$, so that if our tropical virtual dimension is $m$, $\totb{\theta}$  intersects the image of $\totb{ev}_{\bm\alpha}$ transversely restricted to each stratum. This ensures that only curves with general tropical part contribute to the integral $\int_{[\mathcal M_{\bm\alpha}(\ex X)]}\lambda^{2g-2+n}ev_{\bm\alpha}^{*}\theta$, and this integral becomes a sum of contributions for each such general tropical curve with edges restricted to $\totb{\theta}$. Theorem 4.7 of \cite{egw} allows us to write the contribution of $\gamma$ in an equation that looks as follows.
\[k_{\gamma}\int_{\ex X_{\gamma}}\theta\bigwedge_{v\in V(\gamma)}\eta_{v}\] 
Let us explain the terms in this formula a little.  $k_{\gamma}$ is a combinatorial factor as in equation (\ref{kgamma}).   The $\ex X_{\gamma}$ above is the product of a manifold $\ex X_{\gamma_{e}}$ for every edge of $\gamma$, where $\ex X_{\gamma_{e}}$ keeps track of the position of the corresponding edge of a holomorphic curve with tropical part isomorphic to $\gamma$. When $\alpha_{e}=0$, this $\ex X_{\gamma_{e}}$  is the inverse image of the edge under  $\ex X\longrightarrow \totb{\ex X}$, isomorphic to  $\mathbb C^{d}\times(\mathbb C^{*})^{3-d}$, with $d$ depending on the position of $e$ in $\totb{\ex X}$. When $\alpha_{e}\neq 0$, this $\ex X_{\gamma_{e}}$ is the quotient of the inverse image of the center of the edge under $\ex X\longrightarrow \totb{\ex X}$ by the (free and proper) action of $\mathbb C^{*}$ with weight $\alpha_{e}/\abs\alpha_{e}$. Such a $\ex X_{\gamma_{e}}$ is isomorphic to $\mathbb C^{d}\times(\mathbb C^{*})^{2-d}$.  Compatible with keeping track of the position of the relevant edges, there are natural maps 
\[\ex X_{\gamma}\longrightarrow \ex X_{\bm \alpha}\] 
and
\[\ex X_{\gamma}\longrightarrow (\ex X\tc v)_{\bm\alpha_{v}}\ .\]
To explain $\eta_{v}$,  let $\bm\alpha_{v}$ be the list of derivatives of edges leaving $v$, and consider \[ev_{\bm\alpha_{v}}:[\Msw_{\bm\alpha_{v}}(\ex X\tc v)]\longrightarrow (\ex X\tc v)_{\bm\alpha_{v}}\ .\]  
The  $\eta_{v}$ in the integral above is the pull back of  $(ev_{\bm\alpha_{v}})_{!}\lambda^{2g-2+n}$ to $\ex X_{\gamma}$ under the natural map above, and $\theta$ in the integral above really means the pullback to $\ex X_{\gamma}$ of $\theta\in \rh^{*}(\ex X_{\bm\alpha})$. The above integral is over a non-compact space (isomorphic to $\mathbb C^{a}\times(\mathbb C^{*})^{b}$), however the integral is well defined for $\theta\in\rhf^{*}(\ex X_{\bm\alpha})$, and it is possible to choose $\theta$ so that the integrand is compactly supported.

  We can similarly describe the contribution, to Gromov-Witten invariants of $\ex Y$, of the tropical curve $\gamma'$ in $\totb{\ex Y}$ obtained from $\gamma$ by adding edges with special derivative until $\gamma'$ becomes balanced. This contribution looks like
\[k_{\gamma'}\int_{\ex Y_{\gamma'}}\theta\bigwedge_{v\in V(\gamma)}\eta'_{v}\] 
where $k_{\gamma'}\abs{\Aut\gamma'}=k_{\gamma}\abs{\Aut\gamma}$ because $\gamma'$ has the same internal edges as $\gamma$.  The spaces involved in the above expressions are related by the following natural commutative diagram.
\[\begin{tikzcd}\ex Y_{(\bm d,\bm\alpha)}\dar&\lar \ex Y_{\gamma'}\dar\rar &\ex Y_{(\bm d_{v},\bm\alpha_{v})}\dar
\\ \ex X_{\bm\alpha}& \lar \ex X_{\gamma}\rar &\ex X_{\bm\alpha_{v}}\end{tikzcd}\]

 Because our vertices satisfy the correspondence formula (\ref{gac}) and therefore (\ref{eta c}), we can push forward this integral over $\ex Y_{\gamma'}$ to $\ex X_{\gamma}$ to obtain
\[k_{\gamma}\frac{\abs{\Aut\gamma}}{\abs{\Aut\gamma'}}\int_{\ex X_{\gamma}}\theta\bigwedge_{v\in V(\gamma)}(\eta_{v}\prod_{S}\frac{d_{v,S}!}{F^{d_{v,S}}})\]
where $\gamma'$ has $d_{v,S}$ extra edges attached to $v$ with derivative $\alpha_{S}$. The above integral is $\abs{\Aut\gamma}(\prod_{v,S}d_{v,S}!)/F^{\sum_{S} d_{S}}\abs{\Aut\gamma'}$ times the contribution of $\gamma$ to our Gromov-Witten invariant of $\ex X$. There are $\prod_{S}d_{S}!/\prod_{v,S}d_{v,S}!$ distinct ways of distributing labels to the extra edges in $\gamma'$, so the contribution of all these curves to our Gromov-Witten invariant of $\ex Y$ is $\prod_{S}d_{S}!/F^{d_{S}}$ times the contribution of $\gamma$ to our Gromov-Witten invariant of $\ex X$. (In the case that $\Aut\gamma'\neq\Aut\gamma$, the automorphisms of $\gamma'$ are the automorphisms of $\gamma$ that fix each unbalanced vertex. The relationship between the contribution of $\gamma$ and the contributions of the corresponding $\gamma'$ are still valid in this case.) 

To verify that our Gromov-Witten invariant of $\ex Y$ is $\prod_{S}d_{S}!/F^{d_{S}}$ times our Gromov-Witten invariant of $\ex X$, it remains to check that every tropical curve contributing to $\int_{[\mathcal M_{(d,\bm\alpha)}(\ex Y)]}ev_{(\bm d,\bm\alpha)}^{*}\theta$ is in the form of a $\gamma'$ created from a general $\gamma$ in $\totb{\ex X}$. For given generic constraints $\totb\theta$, there are a finite number of tropical curves in $\totb {\ex Y}$ satisfying this constraint and contributing to the above Gromov-Witten invariant. By shrinking this constraint towards $0\in \totb{\ex Y}$,   we may ensure that the only edges intersecting the strata of $\totb{\ex Y}$ not in $\totb{\ex X}$ have special derivative.  These curves are in the form of $\gamma'$ for a general $\gamma$ in $\totb {\ex X}$, so our  Gromov-Witten invariant of $\ex Y$ is $\prod_{S}d_{S}!/F^{d_{S}}$ times our Gromov-Witten invariant of $\ex X$. As $\ex Y$ is a refinement of $\ex T^{3}$, our correspondence formula (\ref{correspondence}) applies to calculate the Gromov-Witten invariants of $\ex Y$, and it follows that formula (\ref{gac}) is satisfied for our Gromov-Witten invariant of $\ex X$. This completes the proof of Claim \ref{induct}.

\

Let us use Claim \ref{induct}  to prove that our correspondence formula (\ref{gac}) holds for $\ex X$ with at most $3$ special vectors, which we shall assume are in the  negative coordinate directions. (This assumption holds for all our  $\ex X\tc v$, after possibly changing coordinates.) Given a collection of ends $\bm\alpha$, construct $3$ polyhedra $P_{\bm\alpha,i}$ as follows: First add special vectors to obtain a  balanced collection of vectors $(\bm d,\bm\alpha)$. Let $P_{\bm\alpha,i}$ be the (non-strictly)  convex polyhedron with edges the projection of vectors in $(\bm d,\bm\alpha)$ that forgets the $i$th coordinate. We shall prove our correspondence formula holds for $\ex X$ by induction on the size of $P_{\bm\alpha,i}$. The smallest cases when $\bm{\alpha}$ does not consist of zero-vectors is when $\bm\alpha$ consists of a single unit vector in one of the coordinate directions. Then one $P_{\bm\alpha,i}$ is $0$ and the other two are unit length intervals. The only general curves with such $\bm\alpha$ have  a  single mono-valent vertex $v$ on a $2$-dimensional stratum of the boundary of $\totb{\ex X}$ and the derivative leaving this vertex is minus the relevant special vector.  Here,  $\ex X\tc v$ is isomorphic to $\ex T^{2}\times\expl (\mathbb CP^{1},0)$, which is the case in which we have verified our gluing formula.  

  Now consider a general tropical curve $\gamma$ in $\totb{\ex X}$ with ends $\bm\alpha$ and a vertex $v$ on the boundary,  and let $\gamma'$ be the corresponding balanced tropical curve. Let $\bm\alpha_{v}$ be the list of derivatives of edges of $\gamma'$ leaving $v$. Then $P_{\bm\alpha_{v},i}$ is contained inside $P_{\bm\alpha,i}$, and strictly contained unless the projection of $\gamma'$ forgetting the ith coordinate is a cone around $v$. (To see this,   consider  the dual graph to an embedded tropical curve with the same image as the projection of $\gamma'$, as in the example pictured  below.) 
 
 \includegraphics{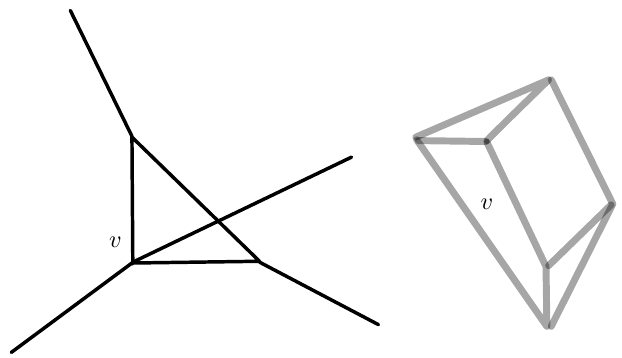}
 
 The only case in which all $P_{\bm\alpha_{v},i}$ are as big as $P_{\bm\alpha,i}$ is when when $\bm\alpha_{v}$ is minimal as considered above, and $\bm\alpha$ is $\bm \alpha_{v}$, possibly with some extra zero-vectors. Claim \ref{induct} gives that our correspondence formula holds for such $\bm\alpha$ because all other vertices must be off the boundary.  By induction on the size of the $P_{\bm\alpha, i }$,  we may therefore assume that our correspondence formula holds at all vertices, so Claim \ref{induct} implies that it also holds for $\bm\alpha$. 
 
 We have now proved our correspondence  formula (\ref{gac}) for all exploded manifolds in the form $\ex X\tc v$.  Claim \ref{induct} therefore implies that our correspondence formula also holds for $\ex X$.

\subsection{Explicit calculation of special cases}

\

\label{zero genus section}

 Consider maps $\mathbb CP^{1}\setminus(0,1,\infty)\longrightarrow (\mathbb C^{*})^{3}$ with data $\alpha$, $\beta$ and $-\alpha-\beta$ at $0,1,\infty$ respectively. Such maps to $(\mathbb C^{*})^{3}$ may be written explicitly as 
\[  (c_{1}z^{\alpha_{1}}(z-1)^{\beta_{1}},c_{2}z^{\alpha_{2}}(z-1)^{\beta_{2}},c_{3}z^{\alpha_{3}}(z-1)^{\beta_{3}})\]
 and the moduli space is naturally parametrized by $(c_{1},c_{2},c_{3})\in (\mathbb C^{*})^{3}$.  The moduli space of maps $\expl(\mathbb CP^{1},\{0,1,\infty\})\longrightarrow \ex T^{3}$ may be written in the same way, except $c_{i}\in \ex T^{3}$, and $z$ and $(z-1)$ now indicate  maps $\expl(\mathbb CP^{1},\{0,1,\infty\})\longrightarrow \ex T$. Moreover, these are the only stable zero genus holomorphic curves in $\ex T^{3}$ with data $\alpha,\beta,(-\alpha-\beta)$. The linearization of the $\dbar$ operator is surjective because it can be identified as the usual $\dbar$ operator acting on $\mathbb C^{3}$-valued functions on $\mathbb CP^{1}$. The moduli space is therefore unobstructed, and  the virtual fundamental class of zero genus curves is parametrized by $\ex T^{3}$ in the above way. Accordingly, $F_{\alpha,\beta,-\alpha-\beta}=\lambda+\dotsc$.
 
 A special case of the above calculation is when $\beta=0\neq\alpha$. In this case, all higher genus invariants are $0$. This is because with fewer than $3$ ends, the virtual dimension of the moduli stack of curves is less than $6$, which is the minimal dimension for curves of positive genus, because the action of $\ex T^{3}$ is free.  It follows that the virtual fundamental class of the moduli stack of such positive genus curves with data $\alpha,-\alpha$ is $0$, and  the divisor equation (which is valid in the context of exploded manifolds) implies that the virtual fundamental class of positive genus curves with data $\alpha,0,-\alpha$ must also be $0$. Accordingly, $F_{\alpha,0,-\alpha}=\lambda$.

\bibliographystyle{plain}
\bibliography{ref}
 \end{document}